\documentclass[10pt]{article}
\usepackage{amsmath,amsthm,amssymb,amsfonts,latexsym,verbatim,bbm,mathrsfs}
\newtheorem{lemma}{Lemma}

\newtheorem{theo}{Theorem}

\begin{document}

\begin{center}
  \huge{\bf Asymptotic analysis for the ratio of the random sum of squares
  to the square of the random sum with applications to risk measures}
\end{center}

\vspace{0.5cm}

\begin{center}
\large{\bf Sophie A. Ladoucette \footnote{Catholic University of
Leuven, Department of Mathematics, W. de Croylaan 54, B-3001
Leuven, Belgium. Email: sophie.ladoucette@wis.kuleuven.be --
Supported by the grant BDB-B/04/03 of the Catholic University of
Leuven} \hspace{1cm} Jef L. Teugels \footnote{Catholic University
of Leuven, Department of Mathematics, W. de Croylaan 54, B-3001
Leuven, Belgium and EURANDOM, Technological University of
Eindhoven, P.O. Box 513, NL-5600 MB Eindhoven, The Netherlands.
Email: jef.teugels@wis.kuleuven.be}}
\end{center}

\vspace{2cm}

\begin{center}
  {\bf 27/10/2005}
\end{center}

\vspace{2cm}

\textbf{Abstract:} Let $\{X_1, X_2, \ldots\}$ be a sequence of
independent and identically distributed positive random variables
of Pareto-type with index $\alpha>0$ and let $\{N(t); \, t\geq
0\}$ be a counting process independent of the $X_i$'s. For any
fixed $t\geq 0$, define:
\begin{equation*}
  T_{N(t)}:=\frac{X_1^2 + X_2^2 + \cdots + X_{N(t)}^2}
  {\left(X_1 + X_2 + \cdots + X_{N(t)}\right)^2}
\end{equation*}
if $N(t)\geq 1$ and $T_{N(t)}:=0$ otherwise.\\

We derive limiting distributions for $T_{N(t)}$ by assuming some
convergence properties for the counting process. This is even
achieved when both the numerator and the denominator defining
$T_{N(t)}$ exhibit an erratic behavior ($\mathbb{E}X_1=\infty$) or
when only the numerator has an erratic behavior
($\mathbb{E}X_1<\infty$ and $\mathbb{E}X_1^2=\infty$). Thanks to
these results, we obtain asymptotic properties pertaining to both
the sample coefficient of variation and the sample dispersion.

\vspace{4cm}


\vspace{0.5cm}


\clearpage

\section{Introduction}

Let $\{X_1, X_2, \ldots\}$ be a sequence of independent and
identically distributed (i.i.d.) positive random variables with
distribution function $F$ and let $\{N(t); \, t\geq 0\}$ be a
counting process independent of the $X_i$'s. For any fixed $t\geq
0$, define the random variable $T_{N(t)}$ by:
\begin{equation}\label{defT}
  T_{N(t)}:=\frac{X_1^2 + X_2^2 + \cdots + X_{N(t)}^2}{\left(X_1 + X_2 + \cdots +
  X_{N(t)}\right)^2}
\end{equation}
if $N(t)\geq 1$ and $T_{N(t)}:=0$ otherwise.\\

The limiting behavior of arbitrary moments of $T_{N(t)}$ is
derived in Ladoucette~\cite{L05} under the conditions that the
distribution function $F$ of $X_1$ is of \textit{Pareto-type} with
positive index $\alpha$ and that the counting process $\{N(t); \,
t\geq 0\}$ is mixed Poisson. In this paper, we focus on
\textit{convergence in distribution} in deriving limits for the
appropriately normalized random variable $T_{N(t)}$. We therefore
generalize earlier results by Albrecher and Teugels~\cite{at04}
where the counting process is non-random (referred below as the
deterministic case). Here, we still assume that $F$ is of
Pareto-type with positive index $\alpha$ except for one result
where the assumption that the fourth moment of $X_1$ exists is
made. Our results are derived under the extra condition that the
counting process $\{N(t); \, t\geq 0\}$ either
\textit{$\mathcal{D}$-averages in time} or \textit{$p$-averages in
time} according to the range of $\alpha$. The appropriate
definitions along with some properties are given in
Section~\ref{sec_preli}.\\

The results of the paper are mostly obtained by using the theory
of functions of \textit{regular variation} (e.g. Bingham et
al.~\cite{BGT}). Recall that a Pareto-type distribution function
$F$ with index $\alpha>0$ is defined by:
\begin{equation}\label{defPT}
  1-F(x)\sim x^{-\alpha}\ell(x)\quad\text{as}\quad x\rightarrow\infty
\end{equation}
for a slowly varying function $\ell$, and therefore has a
regularly varying tail $1-F$ with index $-\alpha<0$.\\

Let $\mu_\beta$ denote the moment of order $\beta>0$ of $X_1$,
i.e.:
\begin{equation*}
  \mu_\beta:=\mathbb{E}X_1^\beta=\beta\int_0^\infty x^{\beta-1}\left(1-F(x)\right)dx\leq\infty.
\end{equation*}

Clearly, both the numerator and the denominator defining
$T_{N(t)}$ exhibit an erratic behavior if $\mu_1=\infty$, whereas
this is the case only for the numerator if $\mu_1<\infty$ and
$\mu_2=\infty$. When $X_1$ (or equivalently $F$) is of Pareto-type
with index $\alpha>0$, it turns out that $\mu_\beta$ is finite if
$\beta<\alpha$ but infinite whenever $\beta>\alpha$. In
particular, $\mu_1<\infty$ if $\alpha>1$ while $\mu_2<\infty$ as
soon as $\alpha>2$. Since the asymptotic behavior of $T_{N(t)}$ is
influenced by the finiteness of $\mu_1$ and/or $\mu_2$, different
limiting distributions consequently show up according to the range
of $\alpha$. This is expressed in our main results given in
Section~\ref{sec_weaklaws}. In Section~\ref{sec_appli}, we use our
results to study the asymptotic behavior of the \textit{sample
coefficient of variation} and the \textit{sample dispersion}, two
risk measures popular in applications.\\

The \textit{coefficient of variation} of a positive random
variable $X$ with distribution function $F$ is defined by:
\begin{equation*}
  \mathrm{CoVar}(X):=\frac{\sqrt{\mathbb{V}X}}{\mathbb{E}X}
\end{equation*}
where $\mathbb{V}X$ denotes the variance of $X$. This risk measure
is frequently used in practice and is of particular interest to
actuaries since it measures the risk associated with $X$. From a
random sample $X_1, \ldots, X_{N(t)}$ from $X$ of random size
$N(t)$ from a nonnegative integer-valued distribution, the
coefficient of variation $\mathrm{CoVar}(X)$ is naturally
estimated by the sample coefficient of variation defined by:
\begin{equation}\label{defScov}
  \widehat{\mathrm{CoVar}(X)}:=\frac{S}{\overline{X}}
\end{equation}
where $\overline{X}:=\frac{1}{N(t)}\sum_{i=1}^{N(t)}X_i$ is the
sample mean and
$S^2:=\frac{1}{N(t)}\sum_{i=1}^{N(t)}\left(X_i-\overline{X}\right)^2$
is the sample variance.\\

For theoretical and practical results pertaining to the
coefficient of variation within the context of (re)insurance, see
Mack~\cite{MAC97}. The properties of the sample coefficient of
variation are usually studied under the tacite assumption of the
finiteness of sufficiently many moments of $X$. However, the
existence of moments of $X$ is not always guaranteed in practical
applications. It is therefore useful to investigate the limiting
behavior of $\widehat{\mathrm{CoVar}(X)}$ also in these cases. It
turns out that this can be achieved by using results on
$T_{N(t)}$. Indeed, the quantity $T_{N(t)}$ appears as a basic
ingredient in the study of the sample coefficient of variation
since the following holds:
\begin{equation}\label{relcov}
  \widehat{\mathrm{CoVar}(X)}=\sqrt{N(t) \, T_{N(t)}-1}.
\end{equation}

In Subsection~\ref{sec_appli_1}, we take advantage from this link
to derive asymptotic properties of the sample coefficient of
variation under the same assumptions on $X$ and on the counting
process $\{N(t); \, t\geq 0\}$ as in Section~\ref{sec_weaklaws}.
Note that this is done even when the first moment and/or the
second moment of $X$ do not exist.\\

Another risk measure of the positive random variable $X$ that is
very popular is the \textit{dispersion} defined by:
\begin{equation*}
  \mathrm{D}(X):=\frac{\mathbb{V}X}{\mathbb{E}X}.
\end{equation*}

In a (re)insurance context, the value of the dispersion is used to
compare the volatility of a portfolio with respect to the Poisson
case for which the dispersion equals $1$. Similarly to the
coefficient of variation, the dispersion $\mathrm{D}(X)$ is
typically estimated by the sample dispersion defined by:
\begin{equation}\label{defSdisp}
  \widehat{\mathrm{D}(X)}:=\frac{S^2}{\overline{X}}.
\end{equation}

Defining the random variable $C_{N(t)}$ for any fixed $t\geq 0$
by:
\begin{equation}\label{defC}
  C_{N(t)}:=\frac{X_1^2 + X_2^2 + \cdots + X_{N(t)}^2}{X_1 + X_2 + \cdots + X_{N(t)}}
\end{equation}
if $N(t)\geq 1$ and $C_{N(t)}:=0$ otherwise, leads to the
following link with the sample dispersion:
\begin{equation}\label{reldisp}
  \widehat{\mathrm{D}(X)}=C_{N(t)}-\overline{X}.
\end{equation}

It turns out that results from Section~\ref{sec_weaklaws} can be
used to derive asymptotic properties of the sample dispersion from
those of $C_{N(t)}$. The results are given in
Subsection~\ref{sec_appli_2} by using the same conditions on $X$
and on the counting process $\{N(t); \, t\geq 0\}$ as in
Section~\ref{sec_weaklaws}. As for the sample coefficient of
variation, cases where the first moments of $X$ do not exist are
also considered.\\

Finally, Section~\ref{ccl} contains a few conclusions.

\section{Preliminaries}\label{sec_preli}

Though standard notations, we mention that
$\stackrel{a.s.}{\longrightarrow}$,
$\stackrel{p}{\longrightarrow}$,
$\stackrel{\mathcal{D}}{\longrightarrow}$ stand for convergence
almost surely, in probability and in distribution, respectively.
Equality in distribution is denoted by
$\stackrel{\mathcal{D}}{=}$. For a measurable function $f$, we
write $f(x)=o(1)$ if $f(x)\rightarrow 0$ as $x\rightarrow\infty$.
Finally, $\Gamma(.)$ denotes the gamma function.\\

Let $\{N(t); \, t\geq 0\}$ be a counting process. For any fixed
$t\geq 0$, the probability generating function $Q_t(.)$ of the
random variable $N(t)$ is defined by:
\begin{equation*}
  Q_t(z):=\mathbb{E}\big\{z^{N(t)}\big\}=\sum_{n=0}^\infty \mathbb{P}[N(t)=n]\,z^n, \quad |z|\leq 1.
\end{equation*}

Most of our results are obtained by assuming that the counting
process satisfies the following condition:
\begin{equation*}
  \frac{N(t)}{t}\stackrel{\mathcal{D}}{\longrightarrow}\Lambda\quad\text{as}\quad t\rightarrow\infty
\end{equation*}
where the limiting random variable $\Lambda$ is such that
$\mathbb{P}[\Lambda>0]=1$. The counting process is then said to
$\mathcal{D}$-average in time to the random variable $\Lambda$. In
a single case, we however need to require the stronger condition
that the above convergence holds in probability rather than in
distribution, i.e.:
\begin{equation*}
  \frac{N(t)}{t}\stackrel{p}{\longrightarrow}\Lambda\quad\text{as}\quad t\rightarrow\infty
\end{equation*}
in which case the counting process is said to $p$-average in time
to the positive random variable $\Lambda$.\\

Whether the counting process $\mathcal{D}$-averages in time or
$p$-averages in time, there results that
$N(t)\stackrel{p}{\longrightarrow}\infty$ as $t\rightarrow\infty$.
Note that the convergence even holds almost surely by
monotonicity.\\

The convergence in distribution is equivalent to the pointwise
convergence of the corresponding Laplace transforms. One therefore
has for a counting process which $\mathcal{D}$-averages in time or
$p$-averages in time to the random variable $\Lambda$ that:
\begin{equation}\label{lapNt}
  \lim_{t\rightarrow\infty}\mathbb{E}\left\{e^{-\theta\frac{N(t)}{t}}\right\}
  =\mathbb{E}\left\{e^{-\theta\Lambda}\right\}, \quad \theta\geq 0.
\end{equation}

For every $\theta\geq 0$, define $u_\theta(x):=e^{-\theta x}$ for
$x\geq 0$. The family of functions
$\left\{u_\theta\right\}_{\theta\geq 0}$ being equicontinuous, it
turns out that the convergence in~(\ref{lapNt}) holds uniformly in
every finite $\theta$-interval (e.g. Corollary page 252 of
Feller~\cite{Feller2}).\\

Very popular counting processes $\mathcal{D}$-average in time. We
quote a number of examples. The deterministic case provides a
first example where $\Lambda$ is degenerate at the point $1$. Any
mixed Poisson process obviously $\mathcal{D}$-averages in time to
its mixing random variable. We refer to the monograph by
Grandell~\cite{gr97} for a very thorough treatment of mixed
Poisson processes and their properties. Also, any infinitely
divisible distribution on the nonnegative integers generates a
process which $\mathcal{D}$-averages in time as long as the
generic step-size $G$ is a discrete random variable on the
positive integers with $\mathbb{E}G<\infty$. In such a case, the
limiting random variable $\Lambda$ is degenerate at the point
$\lambda\,\mathbb{E}G$, where $\lambda>0$ is the intensity of the
homogeneous Poisson process. Finally, any renewal process
generated by a nonnegative distribution with mean
$\mu\in(0,\infty)$ $\mathcal{D}$-averages in time with $\Lambda$
degenerate at the point $1/\mu$.\\

Let us turn to the process $\{X_i; \, i\geq 1\}$ of i.i.d.
positive random variables with distribution function $F$. As
specified above, most of our results are derived under the
condition that the tail of $F$ satisfies~(\ref{defPT}), i.e. that
$1-F$ is regularly varying with negative index $-\alpha$. Recall
that a measurable and ultimately positive function $f$ on
$\mathbb{R}_+$ is regularly varying with index
$\gamma\in\mathbb{R}$ (written $f\in \mathrm{RV}_\gamma$) if for
all $x>0$, $\lim_{t\rightarrow\infty}f(tx)/f(t)=x^\gamma$. When
$\gamma=0$, $f$ is said to be slowly varying. For a textbook
treatment on the theory of functions of regular variation, we
refer to Bingham et al.~\cite{BGT}.\\

It is well-known that the tail condition~(\ref{defPT}) appears as
the essential condition in the maximal domain of attraction
problem of extreme value theory. For a recent treatment, see
Beirlant et al.~\cite{bgst04}. When $\alpha\in(0,2)$, the
condition is also necessary and sufficient for $F$ to belong to
the additive domain of attraction of a non-normal stable law with
exponent $\alpha$ (e.g. Theorem 8.3.1 of Bingham et
al.~\cite{BGT}).\\

Finally, we give a general result that proves to be useful.

\vspace{0.1cm}

\begin{lemma}\label{lem}
  Let $\{Y_n; \, n\geq 1\}$ be a general sequence of random
  variables and $\{M(t); \, t\geq 0\}$ be a process of nonnegative
  integer-valued random variables. Assume that $\{Y_n; \, n\geq 1\}$ and
  $\{M(t); \, t\geq 0\}$ are independent and that
  $M(t)\stackrel{p}{\longrightarrow}\infty$ as $t\rightarrow\infty$. If
  $Y_n\stackrel{\mathcal{D}}{\longrightarrow}Y$ as $n\rightarrow\infty$
  then $Y_{M(t)}\stackrel{\mathcal{D}}{\longrightarrow}Y$ as $t\rightarrow\infty$.
\end{lemma}

\vspace{0cm}

\begin{proof}
To prove the result of the lemma, we show that
$\lim_{t\rightarrow\infty}\mathbb{P}[Y_{M(t)}\leq y]=F_Y(y)$ at
all continuity points $y$ of the distribution function $F_Y$ of
$Y$.\\

Let $y$ be a point of continuity of $F_Y$. For every
$\epsilon\in(0,1)$, there exists
$n_0=n_0(\epsilon,y)\in\mathbb{N}$ such that if $n>n_0$ then
$\left|\mathbb{P}[Y_n\leq y]-F_Y(y)\right|\leq\epsilon$ since
$Y_n\stackrel{\mathcal{D}}{\longrightarrow}Y$ as
$n\rightarrow\infty$. Furthermore, the independence of the two
processes leads to $\mathbb{P}[Y_{M(t)}\leq y]=\sum_{n=0}^\infty
\mathbb{P}[Y_n\leq y]\,\mathbb{P}[M(t)=n]$ so that:
\begin{eqnarray*}
  \left|\mathbb{P}[Y_{M(t)}\leq y]-F_Y(y)\right| &=&
  \left|\left(\sum_{n=0}^{n_0}+\sum_{n=n_0+1}^\infty\right)
  \left\{\mathbb{P}[Y_n\leq y]-F_Y(y)\right\}\mathbb{P}[M(t)=n]\right| \cr
  &\leq& \sum_{n=0}^{n_0}\left|\mathbb{P}[Y_n\leq y]-F_Y(y)\right|\mathbb{P}[M(t)=n]
  +\sum_{n=n_0+1}^\infty\left|\mathbb{P}[Y_n\leq y]-F_Y(y)\right|\mathbb{P}[M(t)=n] \cr
  &\leq& \mathbb{P}[M(t)\leq n_0]+\epsilon\,\mathbb{P}[M(t)>n_0].
\end{eqnarray*}

Since $M(t)\stackrel{p}{\longrightarrow}\infty$ as
$t\rightarrow\infty$, it consequently follows that
$\limsup_{t\rightarrow\infty}\left|\mathbb{P}[Y_{M(t)}\leq
y]-F_Y(y)\right|\leq\epsilon$ which proves the claim upon letting
$\epsilon\downarrow 0$.
\end{proof}

\section{Convergence in Distribution for $T_{N(t)}$}\label{sec_weaklaws}

We provide limiting distributions for the properly normalized
random variable $T_{N(t)}$ defined in~(\ref{defT}) under the
condition that the distribution function $F$ of $X_1$ is of
Pareto-type with index $\alpha>0$. The last result is even
established by assuming that $\mu_4<\infty$ and consequently holds
in the cases $\alpha=4$ with $\mu_4<\infty$ and $\alpha>4$.
Throughout the section, the counting process $\{N(t); \, t\geq
0\}$ is assumed to $\mathcal{D}$-average in time except for one
result where we need to make the stronger assumption that it
$p$-averages in time.

\vspace{0.1cm}

\paragraph{Case 1: $\alpha\in(0,1)$.}

We start with the case $\alpha\in(0,1)$ in the following theorem.

\vspace{0.1cm}

\begin{theo}\label{theo01}
  Assume that $X_1$ is of Pareto-type with index $\alpha\in(0,1)$. Let $\{N(t); \, t\geq 0\}$
  $\mathcal{D}$-average in time to the random variable $\Lambda$. Then:
  \begin{equation*}
    T_{N(t)}\stackrel{\mathcal{D}}{\longrightarrow}\frac{U_\alpha}{V_\alpha^2}
    \quad \text{as} \quad t\rightarrow\infty
  \end{equation*}
  where the random vector $(U_\alpha, V_\alpha)^\prime$ has the Laplace transform:
  \begin{equation}\label{jointtheo01}
    \mathbb{E}\left\{e^{-rU_\alpha-sV_\alpha}\right\}
    =\mathbb{E}\left\{e^{-\delta_\alpha(r,s) \, \Lambda}\right\}, \quad r>0, \, s\geq 0
  \end{equation}
  with $\delta_\alpha(r,s):=2 \, e^{\frac{s^2}{4r}} \int_0^\infty
  e^{-\left(u+\frac{s}{2\sqrt{r}}\right)^2}\left(u+\frac{s}{2\sqrt{r}}\right)
  \left(\frac{u}{\sqrt{r}}\right)^{-\alpha} du$.\\

  In particular, the random variable $U_\alpha$ has the Laplace transform:
  \begin{equation}\label{LTUalpha}
    \mathbb{E}\left\{e^{-rU_\alpha}\right\}
    =\mathbb{E}\left\{e^{-r^{\alpha/2}\Gamma(1-\alpha/2) \, \Lambda}\right\}, \quad r\geq 0
  \end{equation}
  and the random variable $V_\alpha$ has the Laplace transform:
  \begin{equation}\label{LTValpha}
    \mathbb{E}\left\{e^{-sV_\alpha}\right\}
    =\mathbb{E}\left\{e^{-s^\alpha\Gamma(1-\alpha) \, \Lambda}\right\}, \quad s\geq 0.
  \end{equation}
\end{theo}

\vspace{0cm}

\begin{proof}
Let $1-F(x)\underset{x\uparrow\infty}{\sim} x^{-\alpha}\ell(x)$
for some $\ell\in\mathrm{RV}_0$ and $\alpha\in(0,1)$. Define a
sequence $(a_t)_{t>0}$ by
$1-F(a_t)\underset{t\uparrow\infty}{\sim}\frac{1}{t}$, i.e.
$\lim_{t\rightarrow\infty}t \, a_t^{-\alpha}\ell(a_t)=1$. From the
independence of the $X_i$'s, we get:
\begin{eqnarray*}
  \mathbb{E}\left\{e^{-r\frac{1}{a_t^2}\sum_{i=1}^{N(t)}X_i^2 -s\frac{1}{a_t}\sum_{i=1}^{N(t)}X_i}\right\}
  &=& \sum_{n=0}^{\infty} \mathbb{P}[N(t)=n] \, \mathbb{E}\left\{e^{-r\frac{1}{a_t^2}\sum_{i=1}^{n}X_i^2
  -s\frac{1}{a_t}\sum_{i=1}^{n}X_i}\right\} \cr
  &=& \sum_{n=0}^{\infty} \mathbb{P}[N(t)=n]
  \left(\mathbb{E}\left\{e^{-r\frac{1}{a_t^2} X_1^2 -s\frac{1}{a_t} X_1}\right\}\right)^n \cr
  &=& Q_t\left(e^{-\frac{\delta_{\alpha,t}(r,s)}{t}}\right)
\end{eqnarray*}
with $\delta_{\alpha,t}(r,s):=-t\,\log\int_0^\infty
e^{-r\left(\frac{x}{a_t}\right)^2
-s\frac{x}{a_t}}\,dF(x)\in[0,\infty)$.\\

Assume $r>0$ and $s\geq 0$. By virtue of Lebesgue's theorem on
dominated convergence, it is clear that $\int_0^\infty
e^{-r\left(\frac{x}{a_t}\right)^2 -s\frac{x}{a_t}} \,
dF(x)\rightarrow 1$ as $t\rightarrow\infty$. As a consequence, we
obtain:
\begin{eqnarray*}
  \delta_{\alpha,t}(r,s) &\underset{t\uparrow\infty}{\sim}& -t \left(\int_0^\infty
  e^{-r\left(\frac{x}{a_t}\right)^2 -s\frac{x}{a_t}} \, dF(x)-1\right) \cr
  &=& t \int_0^\infty \left(1-e^{-r\left(\frac{x}{a_t}\right)^2 -s\frac{x}{a_t}}\right) dF(x) \cr
  &=& t \int_0^\infty \left(1-F(x)\right) e^{-r\left(\frac{x}{a_t}\right)^2 -s\frac{x}{a_t}}
  \left(\frac{2r}{a_t^2} \, x + \frac{s}{a_t}\right) dx \cr
  &=& 2 \, e^{\frac{s^2}{4r}}\int_{\frac{s}{2\sqrt{r}}}^\infty y \, e^{-y^2} \, t\left\{1-F\left(a_t
  \left(\frac{y}{\sqrt{r}}-\frac{s}{2r}\right)\right)\right\} dy.
\end{eqnarray*}

We now apply Potter's theorem (e.g. Theorem 1.5.6 of Bingham et
al.~\cite{BGT}) together with Lebesgue's theorem on dominated
convergence. Then the change of variables
$y=u+\frac{s}{2\sqrt{r}}$ leads to:
\begin{eqnarray*}
  \lim_{t\rightarrow\infty}\delta_{\alpha,t}(r,s) &=&
  2 \, e^{\frac{s^2}{4r}} \int_{\frac{s}{2\sqrt{r}}}^\infty y \, e^{-y^2}
  \left(\frac{y}{\sqrt{r}}-\frac{s}{2r}\right)^{-\alpha} dy \cr
  &=& 2 \, e^{\frac{s^2}{4r}} \int_0^\infty
  e^{-\left(u+\frac{s}{2\sqrt{r}}\right)^2}\left(u+\frac{s}{2\sqrt{r}}\right)
  \left(\frac{u}{\sqrt{r}}\right)^{-\alpha} du \cr
  &=:& \delta_\alpha(r,s).
\end{eqnarray*}

Define $\varphi_t(\theta):=Q_t\left(e^{-\frac{\theta}{t}}\right)
=\mathbb{E}\left\{e^{-\theta\frac{N(t)}{t}}\right\}$ for
$\theta\geq 0$. From~(\ref{lapNt}), we know that
$\varphi_t(\theta)\rightarrow\mathbb{E}\left\{e^{-\theta\Lambda}\right\}=:\varphi(\theta)$
as $t\rightarrow\infty$. Moreover, we have just proved that
$\delta_{\alpha,t}(r,s)\rightarrow \delta_\alpha(r,s)$ as
$t\rightarrow\infty$, and it is clear that $\delta_\alpha(r,s)\in
(0,\infty)$ for any $\alpha\in (0,1)$. Write the following
triangular inequality:
\begin{equation*}
  \left|\varphi_t(\delta_{\alpha,t}(r,s))-\varphi(\delta_\alpha(r,s))\right|
  \leq \left|\varphi_t(\delta_{\alpha,t}(r,s))-\varphi(\delta_{\alpha,t}(r,s))\right|
  +\left|\varphi(\delta_{\alpha,t}(r,s))-\varphi(\delta_\alpha(r,s))\right|.
\end{equation*}

On the one hand,
$\lim_{t\rightarrow\infty}\left|\varphi(\delta_{\alpha,t}(r,s))-\varphi(\delta_\alpha(r,s))\right|=0$
by continuity of $\varphi$. On the other hand, for $t$ large
enough, there exist reals $a,b$ with $0\leq
a<\delta_\alpha(r,s)<b$ such that $\delta_{\alpha,t}(r,s)\in
[a,b]$. Then,
$\lim_{t\rightarrow\infty}\left|\varphi_t(\delta_{\alpha,t}(r,s))-\varphi(\delta_{\alpha,t}(r,s))\right|=0$
if and only if
$\lim_{t\rightarrow\infty}\sup_{\theta\in[a,b]}\left|\varphi_t(\theta)-\varphi(\theta)\right|=0$.
The last equivalence is true since~(\ref{lapNt}) holds uniformly
in every finite $\theta$-interval. Hence, we obtain
$\varphi_t(\delta_{\alpha,t}(r,s))\rightarrow\varphi(\delta_\alpha(r,s))$
as $t\rightarrow\infty$, that is:
\begin{equation*}
  \lim_{t\rightarrow\infty}\mathbb{E}\left\{e^{-r \frac{1}{a_t^2} \sum_{i=1}^{N(t)}X_i^2
  -s \frac{1}{a_t} \sum_{i=1}^{N(t)}X_i}\right\}
  =\mathbb{E}\left\{e^{-\delta_\alpha(r,s) \, \Lambda}\right\}, \quad r>0, \, s\geq 0
\end{equation*}
or equivalently:
\begin{equation*}
  \left(\frac{1}{a_t^2} \sum_{i=1}^{N(t)}X_i^2, \, \frac{1}{a_t}\sum_{i=1}^{N(t)}X_i\right)^\prime
  \stackrel{\mathcal{D}}{\longrightarrow} (U_\alpha, V_\alpha)^\prime \quad\text{as}\quad t\rightarrow\infty
\end{equation*}
where the random vector $(U_\alpha, V_\alpha)^\prime$ has the
Laplace transform~(\ref{jointtheo01}).\\

Using the continuous mapping theorem (e.g. Corollary 1 page 31 of
Billingsley~\cite{Bill68}), we thus deduce:
\begin{equation*}
  T_{N(t)}=\frac{\frac{1}{a_t^2} \sum_{i=1}^{N(t)}X_i^2}{\left(\frac{1}{a_t}\sum_{i=1}^{N(t)}X_i\right)^2}
  \stackrel{\mathcal{D}}{\longrightarrow}\frac{U_\alpha}{V_\alpha^2}\quad\text{as}\quad t\rightarrow\infty.
\end{equation*}

The marginal distribution of $U_\alpha$, which is the limiting
distribution of $\frac{1}{a_t^2} \sum_{i=1}^{N(t)}X_i^2$, is
determined by its Laplace transform:
\begin{equation*}
  \mathbb{E}\left\{e^{-rU_\alpha}\right\}
  =\mathbb{E}\left\{e^{-r^{\alpha/2}\Gamma(1-\alpha/2) \, \Lambda}\right\}, \quad r\geq 0
\end{equation*}
which is obtained by setting $s=0$ in~(\ref{jointtheo01}).\\

To get the Laplace transform of the marginal distribution of
$V_\alpha$, which is the limiting distribution of $\frac{1}{a_t}
\sum_{i=1}^{N(t)}X_i$, just remake the argument in the proof above
with $r=0$. Instead of $\delta_{\alpha,t}(r,s)$, work with
$\delta_{\alpha,t}(s):=-t\,\log \int_0^\infty
e^{-s\frac{x}{a_t}}\,dF(x)\in[0,\infty)$. Assuming $s>0$, we
obtain:
\begin{equation*}
  \delta_{\alpha,t}(s) \underset{t\uparrow\infty}{\sim}
  \frac{ts}{a_t} \int_0^\infty \left(1-F(x)\right) e^{-s\frac{x}{a_t}} \, dx
  =\int_0^\infty  e^{-y} \, t\left\{1-F\left(a_t\frac{y}{s}\right)\right\} dy.
\end{equation*}

By Potter's theorem and Lebesgue's theorem on dominated
convergence, we then get:
\begin{equation*}
  \lim_{t\rightarrow\infty}\delta_{\alpha,t}(s)
  =\int_0^\infty e^{-y} \left(\frac{y}{s}\right)^{-\alpha} dy=s^\alpha \Gamma(1-\alpha)
\end{equation*}
leading to:
\begin{equation*}
  \mathbb{E}\left\{e^{-sV_\alpha}\right\}
  =\mathbb{E}\left\{e^{-s^\alpha\Gamma(1-\alpha) \, \Lambda}\right\}, \quad s\geq 0.
\end{equation*}
This concludes the proof.
\end{proof}

\vspace{0.1cm}

Note that in the special case where $\Lambda$ is degenerate at the
point $1$, the random variables $U_\alpha$ and $V_\alpha$ are
stable with respective exponent $\alpha/2$ and $\alpha$. In
particular, this holds in discrete time i.e. when the counting
process is non-random.

\vspace{0.1cm}

\paragraph{Case 2: $\alpha=1$, $\mu_1=\infty$.}

Our next result deals with the case $\alpha=1$ and $\mu_1=\infty$.

\vspace{0.1cm}

\begin{theo}\label{theo1}
  Assume that $X_1$ is of Pareto-type with index $\alpha=1$ and that $\mu_1=\infty$. Let
  $\{N(t); \, t\geq 0\}$ $\mathcal{D}$-average in time to the random variable $\Lambda$. Then:
  \begin{equation*}
    \left(\frac{a^\prime_t}{a_t}\right)^2 T_{N(t)}
    \stackrel{\mathcal{D}}{\longrightarrow} \frac{U_1}{\Lambda^2}\quad\text{as}\quad t\rightarrow\infty
  \end{equation*}
  where the random vector $(U_1, \Lambda)^\prime$ has the Laplace transform:
  \begin{equation}\label{jointtheo1}
    \mathbb{E}\left\{e^{-rU_1-s\Lambda}\right\}=\mathbb{E}\left\{e^{-\sqrt{r\pi} \, \Lambda}\right\},
    \quad r>0, \, s\geq 0
  \end{equation}
  and where $(a_t)_{t>0}$ is defined by $\lim_{t\rightarrow\infty}t \,
  a_t^{-1}\ell(a_t)=1$ and $(a^\prime_t)_{t>0}$ is defined
  by $\lim_{t\rightarrow\infty}t \, a^{\prime -1}_t
  \tilde{\ell}(a^\prime_t)=1$ with $\tilde{\ell}(x)=\int_0^x
  \frac{\ell(u)}{u} \, du \in \mathrm{RV}_0$.\\

  In particular, the random variable $U_1$ has the Laplace transform~\eqref{LTUalpha} with $\alpha=1$.
\end{theo}

\vspace{0cm}

\begin{proof}
Let $1-F(x)\underset{x\uparrow\infty}{\sim} x^{-1}\ell(x)$ for
some $\ell\in \mathrm{RV}_0$ and $\mu_1=\infty$. Define a sequence
$(a_t)_{t>0}$ by
$1-F(a_t)\underset{t\uparrow\infty}{\sim}\frac{1}{t}$, i.e.
$\lim_{t\rightarrow\infty}t \, a_t^{-1}\ell(a_t)=1$, and a
sequence $(a^\prime_t)_{t>0}$ by $\lim_{t\rightarrow\infty}t \,
a^{\prime -1}_t\tilde{\ell}(a^\prime_t)=1$ with
$\tilde{\ell}(x)=\int_0^x \frac{\ell(u)}{u} \, du \in
\mathrm{RV}_0$. From the independence of the $X_i$'s, we get:
\begin{equation*}
  \mathbb{E}\left\{e^{-r\frac{1}{a_t^2}\sum_{i=1}^{N(t)}X_i^2 -s\frac{1}{a^\prime_t}
  \sum_{i=1}^{N(t)}X_i}\right\} = Q_t\left(e^{-\frac{\delta_t(r,s)}{t}}\right)
\end{equation*}
with $\delta_t(r,s):=-t \, \log \int_0^\infty
e^{-r\left(\frac{x}{a_t}\right)^2 -s\frac{x}{a^\prime_t}} \, dF(x)
\in [0,\infty)$.\\

Assume $r>0$ and $s\geq 0$. By virtue of Lebesgue's theorem on
dominated convergence, it is clear that $\int_0^\infty
e^{-r\left(\frac{x}{a_t}\right)^2 -s\frac{x}{a^\prime_t}} \,
dF(x)\rightarrow 1$ as $t\rightarrow\infty$. As a consequence, we
obtain:
\begin{eqnarray*}
  \delta_t(r,s) &\underset{t\uparrow\infty}{\sim}&
  t \int_0^\infty \left(1-F(x)\right) e^{-r\left(\frac{x}{a_t}\right)^2 -s\frac{x}{a^\prime_t}}
  \left(\frac{2r}{a_t^2} \, x + \frac{s}{a^\prime_t}\right) dx \cr
  &=& 2 \, e^{\frac{s^2}{4r}\left(\frac{a_t}{a^\prime_t}\right)^2}
  \int_{\frac{s}{2\sqrt{r}}\frac{a_t}{a^\prime_t}}^\infty y \, e^{-y^2} \, t\left\{1-F\left(a_t
  \left(\frac{y}{\sqrt{r}}-\frac{s}{2r}\frac{a_t}{a^\prime_t}\right)\right)\right\} dy.
\end{eqnarray*}

The use of the de Bruyn conjugate together with
$\lim_{x\rightarrow\infty}\frac{\ell(x)}{\tilde{\ell}(x)}=0$ leads
to $\frac{a_t}{a^\prime_t}\rightarrow 0$ as
$t\rightarrow\infty$.\\

The uniform convergence theorem for regularly varying functions
(e.g. Theorem 1.5.2 of Bingham et al.~\cite{BGT}) together with
Potter's theorem and Lebesgue's theorem on dominated convergence
thus gives:
\begin{equation*}
  \lim_{t\rightarrow\infty}\delta_t(r,s)=2\sqrt{r}\int_0^\infty e^{-y^2} \, dy=\sqrt{r\pi}.
\end{equation*}

A similar argument as in the proof of Theorem~\ref{theo01} applied
to the convergence of
$Q_t\left(e^{-\frac{\delta_t(r,s)}{t}}\right)$ yields:
\begin{equation}\label{joint1_inter}
  \lim_{t\rightarrow\infty}
  \mathbb{E}\left\{e^{-r\frac{1}{a_t^2}\sum_{i=1}^{N(t)}X_i^2 -s\frac{1}{a^\prime_t}
  \sum_{i=1}^{N(t)}X_i}\right\}
  =\mathbb{E}\left\{e^{-\sqrt{r\pi} \, \Lambda}\right\}, \quad r>0, \, s\geq 0
\end{equation}
or equivalently:
\begin{equation*}
  \left(\frac{1}{a_t^2} \sum_{i=1}^{N(t)}X_i^2, \, \frac{1}{a^\prime_t} \sum_{i=1}^{N(t)}X_i\right)^\prime
  \stackrel{\mathcal{D}}{\longrightarrow} (Y, Z)^\prime \quad\text{as}\quad t\rightarrow\infty
\end{equation*}
where the Laplace transform of the random vector $(Y, Z)^\prime$
is the term on the right-hand side of~(\ref{joint1_inter}).\\

The marginal distribution of $Y$, which is the limiting
distribution of $\frac{1}{a_t^2} \sum_{i=1}^{N(t)}X_i^2$, is
determined by its Laplace transform:
\begin{equation*}
  \mathbb{E}\left\{e^{-rY}\right\}=\mathbb{E}\left\{e^{-\sqrt{r\pi} \, \Lambda}\right\}, \quad r\geq 0
\end{equation*}
which is obtained by setting $s=0$ in~(\ref{joint1_inter}). Hence,
it turns out that $Y\stackrel{\mathcal{D}}{=}U_1$ where $U_1$ has
the Laplace transform~(\ref{LTUalpha}) with $\alpha=1$.\\

To get the Laplace transform of the marginal distribution of $Z$,
which is the limiting distribution of $\frac{1}{a^\prime_t}
\sum_{i=1}^{N(t)}X_i$, just remake the argument in the proof above
with $r=0$. Instead of $\delta_t(r,s)$, work with
$\delta_t(s):=-t\,\log\int_0^\infty
e^{-s\frac{x}{a^\prime_t}}\,dF(x)\in[0,\infty)$. Assuming $s>0$,
we obtain:
\begin{equation*}
  \delta_t(s) \underset{t\uparrow\infty}{\sim}
  \frac{ts}{a^\prime_t}\int_0^\infty \left(1-F(x)\right) e^{-s\frac{x}{a^\prime_t}} \, dx
  =\frac{ts}{a^\prime_t}\int_0^\infty e^{-y}\int_0^{\frac{a^\prime_t y}{s}}\left(1-F(x)\right) dx \, dy.
\end{equation*}

Since $\tilde{\ell}(x)=\int_0^x \frac{\ell(u)}{u} \, du \sim
\int_0^x \left(1-F(u)\right) du$ as $x\rightarrow\infty$ and using
the uniform convergence theorem for slowly varying functions, we
get:
\begin{equation*}
  \int_0^{\frac{a^\prime_t y}{s}} \left(1-F(x)\right) dx \sim
  \tilde{\ell}\left(\frac{a^\prime_t y}{s}\right)\sim
  \tilde{\ell}\left(a^\prime_t\right) \sim \frac{a^\prime_t}{t}\quad\text{as}\quad t\rightarrow\infty.
\end{equation*}

By Lebesgue's theorem on dominated convergence, we thus get
$\lim_{t\rightarrow\infty}\delta_t(s)=s$. This leads to:
\begin{equation*}
  \mathbb{E}\left\{e^{-sZ}\right\}=\mathbb{E}\left\{e^{-s\Lambda}\right\}, \quad s\geq 0
\end{equation*}
and consequently $Z\stackrel{\mathcal{D}}{=}\Lambda$.\\

Therefore, we have:
\begin{equation*}
  \left(\frac{1}{a_t^2} \sum_{i=1}^{N(t)}X_i^2, \, \frac{1}{a^\prime_t} \sum_{i=1}^{N(t)}X_i\right)^\prime
  \stackrel{\mathcal{D}}{\longrightarrow} (U_1, \Lambda)^\prime \quad\text{as}\quad t\rightarrow\infty
\end{equation*}
where the joint distribution function of $(U_1, \Lambda)^\prime$
is given through~(\ref{jointtheo1}).\\

Using the continuous mapping theorem, we thus deduce:
\begin{equation*}
  \left(\frac{a^\prime_t}{a_t}\right)^2 T_{N(t)}=\frac{\frac{1}{a_t^2}
  \sum_{i=1}^{N(t)}X_i^2}{\left(\frac{1}{a^\prime_t}\sum_{i=1}^{N(t)}X_i\right)^2}
  \stackrel{\mathcal{D}}{\longrightarrow}\frac{U_1}{\Lambda^2} \quad\text{as}\quad
  t\rightarrow\infty
\end{equation*}
and the proof is complete.
\end{proof}

\vspace{0.1cm}

\paragraph{Case 3: $\alpha\in(1,2)$ or $\alpha=1$, $\mu_1<\infty$.}

In the following theorem, the case $\alpha\in(1,2)$ (including
$\alpha=1$ if $\mu_1<\infty$) is treated. Two different results
are proved according to the choice for the normalization of
$T_{N(t)}$.

\vspace{0.1cm}

\begin{theo}\label{theo12}
  Assume that $X_1$ is of Pareto-type with index $\alpha\in(1,2)$ (including $\alpha=1$ if
  $\mu_1<\infty$). Let $\{N(t); \, t\geq 0\}$ $\mathcal{D}$-average in time to the random
  variable $\Lambda$.
  \begin{itemize}
  \item[$(a)$] Then:
  \begin{equation*}
    \left(\frac{N(t)}{a_t}\right)^2 T_{N(t)}\stackrel{\mathcal{D}}{\longrightarrow}
    \frac{1}{\mu_1^2} \, U_\alpha \quad\text{as}\quad t\rightarrow\infty
  \end{equation*}
  where the random variable $U_\alpha$ has the Laplace transform~\eqref{LTUalpha}.
  \item[$(b)$] Then:
  \begin{equation*}
    \left(\frac{t}{a_t}\right)^2 T_{N(t)}\stackrel{\mathcal{D}}{\longrightarrow}
    \frac{1}{\mu_1^2} \, \frac{U_\alpha}{\Lambda^2} \quad\text{as}\quad t\rightarrow\infty
  \end{equation*}
  where the random vector $(U_\alpha, \mu_1\Lambda)^\prime$ has the Laplace transform:
  \begin{equation}\label{jointtheo12}
    \mathbb{E}\left\{e^{-rU_\alpha-s\mu_1\Lambda}\right\}=
    \mathbb{E}\left\{e^{-r^{\alpha/2}\Gamma(1-\alpha/2) \, \Lambda}\right\} , \quad r>0, \, s\geq 0.
  \end{equation}
  In particular, the random variable $U_\alpha$ has the Laplace transform~\eqref{LTUalpha}.
  \end{itemize}
  In $(a)$ and $(b)$, the sequence $(a_t)_{t>0}$ is defined by
  $\lim_{t\rightarrow\infty}t \, a_t^{-\alpha}\ell(a_t)=1$.
\end{theo}

\vspace{0cm}

\begin{proof}
Let $1-F(x)\underset{x\uparrow\infty}{\sim} x^{-\alpha}\ell(x)$
for some $\ell\in \mathrm{RV}_0$ and $\alpha\in(1,2)$ or
$\alpha=1$ if $\mu_1<\infty$. Define a sequence $(a_t)_{t>0}$ by
$1-F(a_t)\underset{t\uparrow\infty}{\sim}\frac{1}{t}$, i.e.
$\lim_{t\rightarrow\infty}t \, a_t^{-\alpha}\ell(a_t)=1$.\\

$(a)$ Since $\mu_1<\infty$ and
$N(t)\stackrel{a.s.}{\longrightarrow}\infty$ as
$t\rightarrow\infty$, it follows by Lemma~\ref{lem} that
$\frac{1}{N(t)}\sum_{i=1}^{N(t)}X_i\stackrel{p}{\longrightarrow}\mu_1$
as $t\rightarrow\infty$.\\

Similar arguments as in the proof of Theorem~\ref{theo01} lead to
$\frac{1}{a_t^2} \sum_{i=1}^{N(t)}X_i^2
\stackrel{\mathcal{D}}{\longrightarrow}U_\alpha$, where the random
variable $U_\alpha$ has the Laplace transform~(\ref{LTUalpha}).\\

Hence, Slutsky's theorem (e.g. Corollary page 97 of
Chung~\cite{Chung01}) and the continuous mapping theorem yield:
\begin{equation*}
  \left(\frac{N(t)}{a_t}\right)^2 T_{N(t)}=
  \frac{\frac{1}{a_t^2} \sum_{i=1}^{N(t)}X_i^2}{\left(\frac{1}{N(t)}\sum_{i=1}^{N(t)}X_i\right)^2}
  \stackrel{\mathcal{D}}{\longrightarrow} \frac{1}{\mu_1^2} \, U_\alpha \quad\text{as}\quad
  t\rightarrow\infty.
\end{equation*}

\vspace{0.2cm}

$(b)$ From the independence of the $X_i$'s, we get:
\begin{equation*}
  \mathbb{E}\left\{e^{-r\frac{1}{a_t^2}\sum_{i=1}^{N(t)}X_i^2 -s\frac{1}{t}\sum_{i=1}^{N(t)}X_i}\right\}
  = Q_t\left(e^{-\frac{\delta_{\alpha,t}(r,s)}{t}}\right)
\end{equation*}
with $\delta_{\alpha,t}(r,s):=-t \, \log \int_0^\infty
e^{-r\left(\frac{x}{a_t}\right)^2 -s\frac{x}{t}} \, dF(x) \in
[0,\infty)$.\\

Assume $r>0$ and $s\geq 0$. By virtue of Lebesgue's theorem on
dominated convergence, it is clear that $\int_0^\infty
e^{-r\left(\frac{x}{a_t}\right)^2 -s\frac{x}{t}} \,
dF(x)\rightarrow 1$ as $t\rightarrow\infty$. As a consequence, we
obtain:
\begin{equation*}
  \delta_{\alpha,t}(r,s) \underset{t\uparrow\infty}{\sim}
  2\,e^{\frac{s^2}{4r}\left(\frac{a_t}{t}\right)^2}
  \int_{\frac{s}{2\sqrt{r}}\frac{a_t}{t}}^\infty y\,e^{-y^2}\,t\left\{1-F\left(a_t
  \left(\frac{y}{\sqrt{r}}-\frac{s}{2r}\frac{a_t}{t}\right)\right)\right\} dy.
\end{equation*}

Now, we prove that $\lim_{t\rightarrow\infty}\frac{a_t}{t}=0$.
When $\alpha=1$, we know that $\ell(x)=o(1)$ since $\mu_1<\infty$.
Consequently, $\frac{a_t}{t}\sim\ell(a_t)\rightarrow 0$ as
$t\rightarrow\infty$. When $\alpha\in (1,2)$, we have
$\frac{a_t}{t}\sim a_t^{1-\alpha}\ell(a_t)\rightarrow 0$ as
$t\rightarrow\infty$ since $1-\alpha\in (-1,0)$.\\

The uniform convergence theorem for regularly varying functions
together with Potter's theorem and Lebesgue's theorem on dominated
convergence thus leads to:
\begin{equation*}
  \lim_{t\rightarrow\infty}\delta_{\alpha,t}(r,s)=
  2 \, r^{\alpha/2}\int_0^\infty y^{1-\alpha} e^{-y^2} \, dy=r^{\alpha/2}\Gamma(1-\alpha/2).
\end{equation*}

A similar argument as in the proof of Theorem~\ref{theo01} applied
to the convergence of
$Q_t\left(e^{-\frac{\delta_{\alpha,t}(r,s)}{t}}\right)$ yields:
\begin{equation}\label{joint12_inter}
  \lim_{t\rightarrow\infty}
  \mathbb{E}\left\{e^{-r\frac{1}{a_t^2}\sum_{i=1}^{N(t)}X_i^2 -s\frac{1}{t}\sum_{i=1}^{N(t)}X_i}\right\}
  =\mathbb{E}\left\{e^{-r^{\alpha/2}\Gamma(1-\alpha/2) \, \Lambda}\right\}, \quad r>0, \, s\geq 0
\end{equation}
or equivalently:
\begin{equation*}
  \left(\frac{1}{a_t^2} \sum_{i=1}^{N(t)}X_i^2, \, \frac{1}{t} \sum_{i=1}^{N(t)}X_i\right)^\prime
  \stackrel{\mathcal{D}}{\longrightarrow} (Y, Z)^\prime \quad\text{as}\quad t\rightarrow\infty
\end{equation*}
where the Laplace transform of the random vector $(Y, Z)^\prime$
is the term on the right-hand side of~(\ref{joint12_inter}).\\

The marginal distribution of $Y$, which is the limiting
distribution of $\frac{1}{a_t^2} \sum_{i=1}^{N(t)}X_i^2$, is
determined by its Laplace transform:
\begin{equation*}
  \mathbb{E}\left\{e^{-rY}\right\}
  =\mathbb{E}\left\{e^{-r^{\alpha/2}\Gamma(1-\alpha/2) \, \Lambda}\right\}, \quad r\geq 0
\end{equation*}
which is obtained by setting $s=0$ in~(\ref{joint12_inter}).
Hence, it turns out that $Y\stackrel{\mathcal{D}}{=}U_\alpha$
where $U_\alpha$ has the Laplace transform~(\ref{LTUalpha}).\\

Since $\mu_1<\infty$ and
$N(t)\stackrel{a.s.}{\longrightarrow}\infty$ as
$t\rightarrow\infty$, it follows by Lemma~\ref{lem} that
$\frac{1}{N(t)}\sum_{i=1}^{N(t)}X_i\stackrel{p}{\longrightarrow}\mu_1$
as $t\rightarrow\infty$. Moreover, we know by assumption that
$\frac{N(t)}{t}\stackrel{\mathcal{D}}{\longrightarrow}\Lambda$ as
$t\rightarrow\infty$. Hence, Slutsky's theorem leads to:
\begin{equation*}
  \frac{1}{t}\sum_{i=1}^{N(t)}X_i=\frac{N(t)}{t} \frac{1}{N(t)}\sum_{i=1}^{N(t)}X_i
  \stackrel{\mathcal{D}}{\longrightarrow}\mu_1\Lambda \quad\text{as}\quad t\rightarrow \infty
\end{equation*}
and consequently $Z\stackrel{\mathcal{D}}{=}\mu_1\Lambda$.\\

Therefore, we have:
\begin{equation*}
  \left(\frac{1}{a_t^2} \sum_{i=1}^{N(t)}X_i^2, \, \frac{1}{t} \sum_{i=1}^{N(t)}X_i\right)^\prime
  \stackrel{\mathcal{D}}{\longrightarrow} (U_\alpha, \mu_1\Lambda)^\prime
  \quad\text{as}\quad t\rightarrow\infty
\end{equation*}
where the joint distribution function of $(U_\alpha,
\mu_1\Lambda)^\prime$ is given through~(\ref{jointtheo12}).\\

Using the continuous mapping theorem, we thus deduce:
\begin{equation*}
  \left(\frac{t}{a_t}\right)^2 T_{N(t)}=\frac{\frac{1}{a_t^2}
  \sum_{i=1}^{N(t)}X_i^2}{\left(\frac{1}{t}\sum_{i=1}^{N(t)}X_i\right)^2}
  \stackrel{\mathcal{D}}{\longrightarrow} \frac{1}{\mu_1^2} \, \frac{U_\alpha}{\Lambda^2}
  \quad\text{as}\quad t\rightarrow\infty
\end{equation*}
and the proof is finished.
\end{proof}

\vspace{0.1cm}

\paragraph{Case 4: $\alpha=2$, $\mu_2=\infty$.}

We pass to the case $\alpha=2$ and $\mu_2=\infty$. As in the
preceding theorem, two different results are derived according to
the way in which $T_{N(t)}$ is normalized.

\vspace{0.1cm}

\begin{theo}\label{theo2}
  Assume that $X_1$ is of Pareto-type with index $\alpha=2$ and
  that $\mu_2=\infty$. Let $\{N(t); \, t\geq 0\}$ $\mathcal{D}$-average in time to the
  random variable $\Lambda$.
  \begin{itemize}
  \item[$(a)$] Then:
  \begin{equation*}
    \left(\frac{N(t)}{a^\prime_t}\right)^2 T_{N(t)}\stackrel{\mathcal{D}}{\longrightarrow}
    \frac{2}{\mu_1^2} \, \Lambda \quad\text{as}\quad t\rightarrow\infty.
  \end{equation*}
  \item[$(b)$] Then:
  \begin{equation*}
    \left(\frac{t}{a^\prime_t}\right)^2 T_{N(t)}\stackrel{\mathcal{D}}{\longrightarrow}
    \frac{2}{\mu_1^2} \, \frac{1}{\Lambda} \quad\text{as}\quad t\rightarrow\infty.
  \end{equation*}
  \end{itemize}
  In $(a)$ and $(b)$, the sequence $(a^\prime_t)_{t>0}$ is defined by
  $\lim_{t\rightarrow\infty}t \, a^{\prime -2}_t\tilde{\ell}(a^\prime_t)=1$
  with $\tilde{\ell}(x)=\int_0^x \frac{\ell(u)}{u} \, du \in \mathrm{RV}_0$.
\end{theo}

\vspace{0cm}

\begin{proof}
Let $1-F(x)\underset{x\uparrow\infty}{\sim} x^{-2}\ell(x)$ for
some $\ell\in \mathrm{RV}_0$ and $\mu_2=\infty$. Define a sequence
$(a^\prime_t)_{t>0}$ by $\lim_{t\rightarrow\infty}t \, a^{\prime
-2}_t \tilde{\ell}(a^\prime_t)=1$ with $\tilde{\ell}(x)=\int_0^x
\frac{\ell(u)}{u} \, du \in \mathrm{RV}_0$.\\

$(a)$ First of all, since $\mu_1<\infty$ and
$N(t)\stackrel{a.s.}{\longrightarrow}\infty$ as
$t\rightarrow\infty$, it follows by Lemma~\ref{lem} that
$\frac{1}{N(t)}\sum_{i=1}^{N(t)}X_i\stackrel{p}{\longrightarrow}\mu_1$
as $t\rightarrow\infty$.\\

From the independence of the $X_i$'s, we get:
\begin{equation*}
  \mathbb{E}\left\{e^{-r\frac{1}{a^{\prime 2}_t}\sum_{i=1}^{N(t)}X_i^2}\right\}
  = Q_t\left(e^{-\frac{\delta_t(r)}{t}}\right)
\end{equation*}
with $\delta_t(r):=-t \, \log \int_0^\infty
e^{-r\left(\frac{x}{a^\prime_t}\right)^2} \, dF(x) \in [0,\infty)$.\\

Assume $r>0$. By Lebesgue's theorem on dominated convergence, it
is clear that $\int_0^\infty
e^{-r\left(\frac{x}{a^\prime_t}\right)^2} dF(x)\rightarrow 1$ as
$t\rightarrow\infty$. As a consequence, we obtain:
\begin{equation*}
  \delta_t(r) \underset{t\uparrow\infty}{\sim}
  \frac{2tr}{a^{\prime 2}_t} \int_0^\infty x\left(1-F(x)\right) e^{-r\left(\frac{x}{a^\prime_t}\right)^2} dx
  =\frac{2tr}{a^{\prime 2}_t} \int_0^\infty e^{-y}\int_0^{a^\prime_t\sqrt{\frac{y}{r}}}
  x\left(1-F(x)\right) dx \, dy.
\end{equation*}

Since $\tilde{\ell}(x)=\int_0^x \frac{\ell(u)}{u} \, du \sim
\int_0^x u \left(1-F(u)\right) du$ as $x\rightarrow\infty$, the
uniform convergence theorem for slowly varying functions leads to:
\begin{equation*}
  \int_0^{a^\prime_t\sqrt{\frac{y}{r}}} x \left(1-F(x)\right) dx \sim
  \tilde{\ell}\left(a^\prime_t\sqrt{\frac{y}{r}}\right)\sim
  \tilde{\ell}\left(a^\prime_t\right) \sim \frac{a^{\prime 2}_t}{t} \quad\text{as}\quad t\rightarrow\infty.
\end{equation*}

By Lebesgue's theorem on dominated convergence, we thus get
$\lim_{t\rightarrow\infty}\delta_t(r)=2r$.\\

A similar argument as in the proof of Theorem~\ref{theo01} applied
to the convergence of $Q_t\left(e^{-\frac{\delta_t(r)}{t}}\right)$
yields:
\begin{equation*}
  \lim_{t\rightarrow\infty}\mathbb{E}\left\{e^{-r \frac{1}{a^{\prime 2}_t} \sum_{i=1}^{N(t)}X_i^2}\right\}
  =\mathbb{E}\left\{e^{-2r\Lambda}\right\}, \quad r\geq 0
\end{equation*}
or equivalently:
\begin{equation*}
  \frac{1}{a^{\prime 2}_t} \sum_{i=1}^{N(t)}X_i^2 \stackrel{\mathcal{D}}{\longrightarrow}
  2\Lambda \quad\text{as}\quad t\rightarrow\infty.
\end{equation*}

Slutsky's theorem together with the continuous mapping theorem
finally gives:
\begin{equation*}
  \left(\frac{N(t)}{a^\prime_t}\right)^2 T_{N(t)}
   =\frac{\frac{1}{a^{\prime 2}_t} \sum_{i=1}^{N(t)}X_i^2}{\left(\frac{1}{N(t)}
   \sum_{i=1}^{N(t)}X_i\right)^2} \stackrel{\mathcal{D}}{\longrightarrow}
  \frac{2}{\mu_1^2} \, \Lambda \quad\text{as}\quad t\rightarrow\infty.
\end{equation*}

\vspace{0.2cm}

$(b)$ From the independence of the $X_i$'s, we get:
\begin{equation*}
  \mathbb{E}\left\{e^{-r\frac{1}{a^{\prime 2}_t}\sum_{i=1}^{N(t)}X_i^2
  -s\frac{1}{t}\sum_{i=1}^{N(t)}X_i}\right\} = Q_t\left(e^{-\frac{\delta_{\alpha,t}(r,s)}{t}}\right)
\end{equation*}
with $\delta_{\alpha,t}(r,s):=-t \, \log \int_0^\infty
e^{-r\left(\frac{x}{a^\prime_t}\right)^2 -s\frac{x}{t}} \, dF(x)
\in [0,\infty)$.\\

Assume $r>0$ and $s\geq 0$. By virtue of Lebesgue's theorem on
dominated convergence, it is clear that $\int_0^\infty
e^{-r\left(\frac{x}{a^\prime_t}\right)^2 -s\frac{x}{t}} \,
dF(x)\rightarrow 1$ as $t\rightarrow\infty$. As a consequence, we
obtain:
\begin{eqnarray*}
  \delta_{\alpha,t}(r,s) &\underset{t\uparrow\infty}{\sim}&
  t\int_0^\infty \left(1-F(x)\right) e^{-r\left(\frac{x}{a^\prime_t}\right)^2 -s\frac{x}{t}}
  \left(\frac{2r}{a^{\prime 2}_t} \, x+\frac{s}{t}\right) dx \cr
  &=& \frac{2tr}{a^{\prime 2}_t}\int_0^\infty x \left(1-F(x)\right)
  \int_{r\left(\frac{x}{a^\prime_t}\right)^2 + s\frac{x}{t}}^\infty e^{-y} \, dy \, dx
  + s\int_0^\infty \left(1-F(x)\right) \int_{r\left(\frac{x}{a^\prime_t}\right)^2 + s\frac{x}{t}}^\infty
  e^{-y} \, dy \, dx \cr
  &=& \frac{2tr}{a^{\prime 2}_t} \int_0^\infty e^{-y}\int_0^{y_*} x \left(1-F(x)\right) dx \, dy
  + s\int_0^\infty e^{-y}\int_0^{y_*} \left(1-F(x)\right) dx \, dy
\end{eqnarray*}
where
$y_*:=a^\prime_t\left(\sqrt{\frac{y}{r}+\left(\frac{s}{2r}\frac{a^\prime_t}{t}\right)^2}
-\frac{s}{2r}\frac{a^\prime_t}{t}\right)$.\\

Since $\lim_{t\rightarrow\infty}t \, a^{\prime -2}_t \,
\tilde{\ell}(a^\prime_t)=1$, we have
$\frac{a^\prime_t}{t}\sim\frac{\tilde{\ell}(a^\prime_t)}{a^\prime_t}\rightarrow
0$ as $t\rightarrow\infty$. Consequently, we get $y_*\sim
a^\prime_t\sqrt{\frac{y}{r}}$ as $t\rightarrow\infty$. Since
$\tilde{\ell}(x)=\int_0^x \frac{\ell(u)}{u} \, du \sim \int_0^x u
\left(1-F(u)\right) du$ as $x\rightarrow\infty$, it follows by the
uniform convergence theorem for slowly varying functions that:
\begin{equation*}
  \int_0^{y_*} x \left(1-F(x)\right) dx \sim \tilde{\ell}\left(y_*\right)
  \sim \tilde{\ell}\left(a^\prime_t\sqrt{\frac{y}{r}}\right)
  \sim \tilde{\ell}\left(a^\prime_t\right)
  \sim \frac{a^{\prime 2}_t}{t} \quad\text{as}\quad t\rightarrow\infty.
\end{equation*}

Moreover, since $\mu_1<\infty$, Lebesgue's theorem on dominated
convergence gives $\int_0^{y_*} \left(1-F(x)\right) dx\rightarrow
\mu_1$ as $t\rightarrow\infty$.\\

Consequently, another application of Lebesgue's theorem on
dominated convergence leads to:
\begin{equation*}
  \lim_{t\rightarrow\infty}\delta_{\alpha,t}(r,s)=2r + \mu_1 s.
\end{equation*}

A similar argument as in the proof of Theorem~\ref{theo01} applied
to the convergence of
$Q_t\left(e^{-\frac{\delta_{\alpha,t}(r,s)}{t}}\right)$ yields:
\begin{equation*}
  \lim_{t\rightarrow\infty}
  \mathbb{E}\left\{e^{-r\frac{1}{a^{\prime 2}_t}\sum_{i=1}^{N(t)}X_i^2
  -s\frac{1}{t}\sum_{i=1}^{N(t)}X_i}\right\}
  =\mathbb{E}\left\{e^{-(2r + \mu_1 s) \, \Lambda}\right\}, \quad r>0, \, s\geq 0
\end{equation*}
or equivalently:
\begin{equation*}
  \left(\frac{1}{a^{\prime 2}_t}\sum_{i=1}^{N(t)}X_i^2, \, \frac{1}{t}\sum_{i=1}^{N(t)}X_i\right)^\prime
  \stackrel{\mathcal{D}}{\longrightarrow} (2\Lambda, \mu_1\Lambda)^\prime
  \quad\text{as}\quad t\rightarrow\infty.
\end{equation*}

Using the continuous mapping theorem, we thus deduce:
\begin{equation*}
  \left(\frac{t}{a^\prime_t}\right)^2 T_{N(t)}
   =\frac{\frac{1}{a^{\prime 2}_t} \sum_{i=1}^{N(t)}X_i^2}{\left(\frac{1}{t}
   \sum_{i=1}^{N(t)}X_i\right)^2} \stackrel{\mathcal{D}}{\longrightarrow}
   \frac{2}{\mu_1^2} \, \frac{1}{\Lambda} \quad\text{as}\quad t\rightarrow\infty
\end{equation*}
and this ends the proof.
\end{proof}

\vspace{0.1cm}

\paragraph{Case 5: $\alpha\in(2,4)$ or $\alpha=2$, $\mu_2<\infty$.}

The following theorem concerns the case $\alpha\in(2,4)$
(including $\alpha=2$ if $\mu_2<\infty$) where it is assumed that
the counting process $p$-averages in time.

\vspace{0.1cm}

\begin{theo}\label{theo24}
  Assume that $X_1$ is of Pareto-type with index $\alpha\in(2,4)$ (including $\alpha=2$
  if $\mu_2<\infty$). Let $\{N(t); \, t\geq 0\}$ $p$-average in time to the random
  variable $\Lambda$. Then:
  \begin{equation*}
    \frac{t^{1-2/\alpha}}{\ell^*(t)}\left(N(t) \, T_{N(t)}-\frac{\mu_2}{\mu_1^2}\right)
    \stackrel{\mathcal{D}}{\longrightarrow} \frac{1}{\mu_1^2} \,
    \frac{W_{\frac{\alpha}{2}}}{\Lambda^{1-2/\alpha}} \quad\text{as}\quad t\rightarrow\infty
  \end{equation*}
  where $W_{\frac{\alpha}{2}}$ is a stable random variable with exponent $\alpha/2$
  independent of $\Lambda$ and $\ell^*\in\mathrm{RV}_0$ is given by
  $\ell^*(t)\sim c_t\,t^{-2/\alpha}$ as $t\rightarrow\infty$ with $(c_t)_{t>0}$
  a sequence defined by $\lim_{t\rightarrow\infty}t\,c_t^{-\alpha/2}\ell(\sqrt{c_t})=1$.
\end{theo}

\vspace{0cm}

\begin{proof}
Let $1-F(x)\underset{x\uparrow\infty}{\sim} x^{-\alpha}\ell(x)$
for some $\ell\in\mathrm{RV}_0$ and $\alpha\in(2,4)$ or $\alpha=2$
if $\mu_2<\infty$. For a sequence $(b_t)_{t>0}$ to be defined,
consider the following identity:
\begin{equation*}
  b_t\left(N(t) \, T_{N(t)}-\frac{\mu_2}{\mu_1^2}\right)
  =\underbrace{b_t\left(\frac{N(t)}{X_1+\cdots+X_{N(t)}}\right)^2
  \left(\frac{1}{N(t)}\sum_{i=1}^{N(t)}X_i^2-\mu_2\right)}_{=:A_t}
  +\underbrace{b_t \, \mu_2\left(\frac{N(t)}{X_1+\cdots+X_{N(t)}}\right)^2-\frac{b_t \,
  \mu_2}{\mu_1^2}}_{=:B_t}.
\end{equation*}

Since
$\mathbb{P}\left[X_1^2>x\right]=\mathbb{P}\left[X_1>\sqrt{x}\right]
\underset{x\uparrow\infty}{\sim}x^{-\alpha/2}\ell\left(\sqrt{x}\right)$,
the tail of $X_1^2$ is regularly varying with index $-\alpha/2$.
Hence, we get $c_n^{-1}\left(\sum_{i=1}^n
X_i^2-n\mu_2\right)\stackrel{\mathcal{D}}{\longrightarrow}W_{\frac{\alpha}{2}}$
as $n\rightarrow\infty$ for a stable random variable
$W_{\frac{\alpha}{2}}$ with exponent $\alpha/2$ and a sequence of
normalizing constants $(c_n)_{n\geq 1}$ defined by
$1-F(\sqrt{c_n})\underset{n\uparrow\infty}{\sim}\frac{1}{n}$, i.e.
$c_n\underset{n\uparrow\infty}{\sim}n^{2/\alpha}\ell^*(n)$ with
$\ell^*(n):=\ell^{2/\alpha}\left(\sqrt{c_n}\right)\in\mathrm{RV}_0$.
Since $N(t)\stackrel{a.s.}{\longrightarrow}\infty$ as
$t\rightarrow\infty$, it follows by Lemma~\ref{lem} that:
\begin{equation}\label{AnsCLTW}
  \frac{N(t)^{1-2/\alpha}}{\ell^*(N(t))}\left(\frac{1}{N(t)}\sum_{i=1}^{N(t)}
  X_i^2-\mu_2\right)\stackrel{\mathcal{D}}{\longrightarrow}W_{\frac{\alpha}{2}}
  \quad\text{as}\quad t\rightarrow\infty.
\end{equation}

Similarly, since $n^{-1/2}\left(\sum_{i=1}^n
X_i-n\mu_1\right)\stackrel{\mathcal{D}}{\longrightarrow}\mathrm{N}(0,\sigma^2)$
as $n\rightarrow\infty$ where the random variable
$\mathrm{N}(0,\sigma^2)$ has the normal distribution with mean $0$
and variance $\sigma^2:=\mathbb{V}X_1<\infty$, we also get:
\begin{equation}\label{AnsCLTN}
  \sqrt{N(t)}\left(\frac{1}{N(t)}\sum_{i=1}^{N(t)}
  X_i-\mu_1\right)\stackrel{\mathcal{D}}{\longrightarrow}\mathrm{N}(0,\sigma^2)
  \quad\text{as}\quad t\rightarrow\infty.
\end{equation}

Define the sequence $(b_t)_{t>0}$ by
$b_t:=\frac{t^{1-2/\alpha}}{\ell^*(t)}$. This leads to:
\begin{equation*}
  A_t=\frac{\ell^*(N(t))}{\ell^*(t)}\left(\frac{t}{N(t)}\right)^{1-2/\alpha}
  \left(\frac{N(t)}{X_1+\cdots+X_{N(t)}}\right)^2\frac{N(t)^{1-2/\alpha}}{\ell^*(N(t))}
  \left(\frac{1}{N(t)}\sum_{i=1}^{N(t)} X_i^2-\mu_2\right).
\end{equation*}

Since $\ell^*\in\mathrm{RV}_0$ and
$\frac{N(t)}{t}\stackrel{p}{\longrightarrow}\Lambda$ as
$t\rightarrow\infty$ with $\mathbb{P}[\Lambda>0]=1$, we get
$\frac{\ell^*(N(t))}{\ell^*(t)}\stackrel{p}{\longrightarrow}1$ as
$t\rightarrow\infty$ by using the uniform convergence theorem for
slowly varying functions and the subsequence principle. Since
$\mu_1<\infty$ and $N(t)\stackrel{a.s.}{\longrightarrow}\infty$ as
$t\rightarrow\infty$, it follows by Lemma~\ref{lem} that
$\frac{1}{N(t)}\sum_{i=1}^{N(t)}X_i\stackrel{p}{\longrightarrow}\mu_1$
as $t\rightarrow\infty$. Recalling~(\ref{AnsCLTW}), Slutsky's
theorem together with the continuous mapping theorem therefore
implies that:
\begin{equation*}
  A_t\stackrel{\mathcal{D}}{\longrightarrow}\frac{1}{\mu_1^2}\,
  \frac{W_{\frac{\alpha}{2}}}{\Lambda^{1-2/\alpha}} \quad \text{as} \quad t\rightarrow\infty
\end{equation*}
thanks to the independence of $W_{\frac{\alpha}{2}}$ and
$\Lambda$. To prove the latter claim, we condition on $N(t)$, use
that $\frac{N(t)}{t}$ is $\sigma(N(t))$-measurable and that
$\{N(t); \, t\geq 0\}$ and $\{X_i; \, i\geq 1\}$ are independent
and apply~(\ref{AnsCLTW}) to get:
\begin{eqnarray*}
  & & Y_t:=\mathbb{P}\left[\left.\frac{N(t)}{t}\leq x,
  \frac{N(t)^{1-2/\alpha}}{\ell^*(N(t))}\left(\frac{1}{N(t)}\sum_{i=1}^{N(t)}
  X_i^2-\mu_2\right)\leq y \right\vert N(t)\right] \cr
  & & \phantom{Y_t}
  =\mathbb{E}\Bigg[\Bigg.\mathbbm{1}_{\left\{\frac{N(t)}{t}\leq x\right\}}
  \,\mathbbm{1}_{\left\{\frac{N(t)^{1-2/\alpha}}{\ell^*(N(t))}\left(\frac{1}{N(t)}\sum_{i=1}^{N(t)}
  X_i^2-\mu_2\right)\leq y\right\}}\Bigg\vert \, N(t)\Bigg] \cr
  & & \phantom{Y_t}
  =\mathbbm{1}_{\left\{\frac{N(t)}{t}\leq x\right\}}
  \,\mathbb{P}\left[\left.\frac{N(t)^{1-2/\alpha}}{\ell^*(N(t))}\left(\frac{1}{N(t)}\sum_{i=1}^{N(t)}
  X_i^2-\mu_2\right)\leq y \right\vert N(t)\right] \cr
  & & \phantom{Y_t}
  =\mathbbm{1}_{\left\{\frac{N(t)}{t}\leq x\right\}}
  \,\mathbb{P}\left[\frac{N(t)^{1-2/\alpha}}{\ell^*(N(t))}\left(\frac{1}{N(t)}\sum_{i=1}^{N(t)}
  X_i^2-\mu_2\right)\leq y\right] \cr
  & & \phantom{Y_t}
  \stackrel{\mathcal{D}}{\longrightarrow}
  \mathbbm{1}_{\left\{\Lambda\leq x\right\}}\,\mathbb{P}\left[W_{\frac{\alpha}{2}}\leq y\right]
  \quad\text{as}\quad t\rightarrow\infty
\end{eqnarray*}
at any continuity points $x$ of the distribution function of
$\Lambda$ and $y$ of that of $W_{\frac{\alpha}{2}}$. The sequence
of random variables $\{Y_t; \, t>0\}$ being uniformly integrable,
we use Theorem 5.4 of Billingsley~\cite{Bill68} to obtain:
\begin{equation*}
  \mathbb{P}\left[\frac{N(t)}{t}\leq x, \frac{N(t)^{1-2/\alpha}}{\ell^*(N(t))}
  \left(\frac{1}{N(t)}\sum_{i=1}^{N(t)} X_i^2-\mu_2\right)\leq y\right]
  \rightarrow\mathbb{P}\left[\Lambda\leq x\right]\,\mathbb{P}\left[W_{\frac{\alpha}{2}}\leq y\right]
  \quad\text{as}\quad t\rightarrow\infty
\end{equation*}
and this proves the claim.\\

At the same time, we also have:
\begin{equation*}
  B_t=-\frac{\mu_2}{\mu_1}\left(\frac{N(t)}{X_1+\cdots+X_{N(t)}}\right)^2
  \sqrt{\frac{t}{N(t)}}\frac{t^{1/2-2/\alpha}}{\ell^*(t)}
  \sqrt{N(t)}\left(\frac{1}{N(t)}\sum_{i=1}^{N(t)} X_i-\mu_1\right)
  \left(1 + \frac{1}{\mu_1} \, \frac{1}{N(t)}\sum_{i=1}^{N(t)} X_i\right).
\end{equation*}

Thanks to~(\ref{AnsCLTN}) and since
$\lim_{t\rightarrow\infty}\frac{t^{1/2-2/\alpha}}{\ell^*(t)}=0$,
Slutsky's theorem, the continuous mapping theorem and the
independence of $\mathrm{N}(0,\sigma^2)$ and $\Lambda$ yield
$B_t\stackrel{p}{\longrightarrow}0$ as $t\rightarrow\infty$. Note
that the independence of $\mathrm{N}(0,\sigma^2)$ and $\Lambda$ is
easily proved using the same kind of arguments as above for the
independence of $W_{\frac{\alpha}{2}}$ and $\Lambda$.\\

Hence, we get by another application of Slutsky's theorem:
\begin{equation*}
  \frac{t^{1-2/\alpha}}{\ell^*(t)}\left(N(t) \, T_{N(t)}-\frac{\mu_2}{\mu_1^2}\right)
  =A_t+B_t\stackrel{\mathcal{D}}{\longrightarrow}\frac{1}{\mu_1^2}
  \, \frac{W_{\frac{\alpha}{2}}}{\Lambda^{1-2/\alpha}}
  \quad\text{as}\quad t\rightarrow\infty
\end{equation*}
and the proof is complete.
\end{proof}

\vspace{0.1cm}

\paragraph{Case 6: $\mu_4<\infty$.}

Finally, the case $\mu_4<\infty$ is under consideration in the
next result which hence holds in the cases $\alpha=4$ with finite
fourth moment and $\alpha>4$.

\vspace{0.1cm}

\begin{theo}\label{theo4}
  Assume that $X_1$ is such that $\mu_4<\infty$. Let $\{N(t); \, t\geq 0\}$
  $\mathcal{D}$-average in time to the random variable $\Lambda$. Then:
  \begin{equation*}
    \sqrt{t}\left(N(t) \, T_{N(t)}-\frac{\mu_2}{\mu_1^2}\right)
    \stackrel{\mathcal{D}}{\longrightarrow} \frac{\mathrm{N}(0,\sigma_*^2)}{\sqrt{\Lambda}}
    \quad\text{as}\quad t\rightarrow\infty
  \end{equation*}
  where the random variables $\mathrm{N}(0,\sigma_*^2)$ and $\Lambda$
  are independent, $\mathrm{N}(0,\sigma_*^2)$ having the normal distribution with mean
  $0$ and variance $\sigma_*^2$ defined by:
  \begin{equation}\label{defsigma*}
    \sigma_*^2:=\frac{\mu_4}{\mu_1^4}-\left(\frac{\mu_2}{\mu_1^2}\right)^2
    +4\left(\frac{\mu_2}{\mu_1^2}\right)^3-\frac{4\mu_2\mu_3}{\mu_1^5}.
  \end{equation}
\end{theo}

\vspace{0cm}

\begin{proof}
Let the distribution function $F$ of $X_1$ be such that
$\mu_4<\infty$. From the bivariate Lindeberg-L\'evy central limit
theorem (e.g. Theorem 1.9.1B of Serfling~\cite{SERF}), one deduces
that:
\begin{equation*}
  \sqrt{n}\left(\frac{1}{n}\sum_{i=1}^n\boldsymbol{Y_n}-\boldsymbol{\mu}\right)
  \stackrel{\mathcal{D}}{\longrightarrow} \mathrm{N}(\boldsymbol{0},\boldsymbol{\Sigma})
  \quad\text{as}\quad n\rightarrow\infty
\end{equation*}
where $\boldsymbol{Y_n}=\left(X_i,X_i^2\right)^\prime$,
$\boldsymbol{\mu}=(\mu_1,\mu_2)^\prime$ and
$\mathrm{N}(\boldsymbol{0},\boldsymbol{\Sigma})$ has the bivariate
normal distribution with mean vector $\boldsymbol{0}=(0,0)^\prime$
and covariance matrix $\boldsymbol{\Sigma}$ defined by:
\begin{equation*}
  \boldsymbol{\Sigma}:=\left(
  \begin{array}{cc}
    \mu_2-\mu_1^2 & \mu_3-\mu_1\mu_2 \\
    \mu_3-\mu_1\mu_2 & \mu_4-\mu_2^2 \\
  \end{array}\right).
\end{equation*}

Following the notation in Serfling~\cite{SERF}, we write this as
$\frac{1}{n}\sum_{i=1}^n\boldsymbol{Y_n}$ is
AN$\left(\boldsymbol{\mu},n^{-1}\boldsymbol{\Sigma}\right)$.\\

By the multivariate delta method, the asymptotic normality carries
over to any function
$g\left(\frac{1}{n}\sum_{i=1}^n\boldsymbol{Y_n}\right)=g\left(\frac{1}{n}\sum_{i=1}^n
X_i,\frac{1}{n}\sum_{i=1}^n X_i^2\right)$ where
$g:\mathbb{R}_+\times\mathbb{R}_+\rightarrow\mathbb{R}$ is
continuously differentiable in a neighborhood of
$\boldsymbol{\mu}$, so that
$g\left(\frac{1}{n}\sum_{i=1}^n\boldsymbol{Y_n}\right)$ is
AN$\left(g(\boldsymbol{\mu}),n^{-1}\boldsymbol{J\Sigma
J^\prime}\right)$ with $\boldsymbol{J}=\left(\frac{\partial
g}{\partial x}(\boldsymbol{\mu}),\frac{\partial g}{\partial
y}(\boldsymbol{\mu})\right)$.\\

With the choice $g(x,y)=\frac{y}{x^2}$, we find that $n\,T_n$ is
AN$\left(\frac{\mu_2}{\mu_1^2},\frac{\sigma_*^2}{n}\right)$ with
$\sigma_*^2$ given by~(\ref{defsigma*}). Since
$N(t)\stackrel{a.s.}{\longrightarrow}\infty$ as
$t\rightarrow\infty$, it consequently follows by Lemma~\ref{lem}
that:
\begin{equation*}
  \sqrt{N(t)}\left(N(t)\,T_{N(t)}-\frac{\mu_2}{\mu_1^2}\right)
  \stackrel{\mathcal{D}}{\longrightarrow} \mathrm{N}(0,\sigma_*^2) \quad\text{as}\quad t\rightarrow\infty
\end{equation*}
where the random variable $\mathrm{N}(0,\sigma_*^2)$ has the
normal distribution with mean $0$ and variance $\sigma_*^2$.\\

The continuous mapping theorem together with the independence of
$\mathrm{N}(0,\sigma_*^2)$ and $\Lambda$ (which is proved using
the same arguments as for the independence of
$W_{\frac{\alpha}{2}}$ and $\Lambda$ in the proof of
Theorem~\ref{theo24}) finally gives:
\begin{equation*}
  \sqrt{t}\left(N(t) \, T_{N(t)}-\frac{\mu_2}{\mu_1^2}\right)
  =\sqrt{\frac{t}{N(t)}}\,\sqrt{N(t)}\left(N(t) \, T_{N(t)}-\frac{\mu_2}{\mu_1^2}\right)
  \stackrel{\mathcal{D}}{\longrightarrow}\frac{\mathrm{N}(0,\sigma_*^2)}{\sqrt{\Lambda}}
  \quad\text{as}\quad t\rightarrow\infty.
\end{equation*}
This completes the proof.
\end{proof}

\vspace{0.1cm}

We end by a remark. There is a slight incompleteness in the
results for the case $\alpha=4$. When $\mu_4<\infty$, the result
is given by Theorem~\ref{theo4}. However when $\mu_4=\infty$, the
results depend on the behavior of the slowly varying function
$\ell$ and the arguments get even more complicated, except for
$\ell(x)\rightarrow\infty$ as $t\rightarrow\infty$ since
Theorem~\ref{theo24} holds with a normally distributed random
variable $W_2$. We hope to treat this remaining case in the
future.

\section{Applications to Risk Measures}\label{sec_appli}

Assume that $X$ is a positive random variable with distribution
function $F$ and let $X_1, \ldots, X_{N(t)}$ be a random sample
from $F$ of random size $N(t)$ from a nonnegative integer-valued
distribution. Thanks to the limiting results derived in
Section~\ref{sec_weaklaws} and the relations~(\ref{relcov})
and~(\ref{reldisp}), we investigate the asymptotic behavior of two
popular risk measures through their distributions.
Subsection~\ref{sec_appli_1} deals with the sample coefficient of
variation $\widehat{\mathrm{CoVar}(X)}$ defined in~(\ref{defScov})
and Subsection~\ref{sec_appli_2} concerns the sample dispersion
$\widehat{\mathrm{D}(X)}$ defined in~(\ref{defSdisp}). The results
are obtained under the same assumptions on $X$ and on the counting
process $\{N(t); \, t\geq 0\}$ as in Section~\ref{sec_weaklaws}.

\subsection{Sample Coefficient of Variation}\label{sec_appli_1}

We determine limits in distribution for the appropriately
normalized random variable $\widehat{\mathrm{CoVar}(X)}$ by using
the distributional results derived in Section~\ref{sec_weaklaws}
for $T_{N(t)}$ and thanks to~(\ref{relcov}). Consequently,
different cases arise depending on the range of $\alpha$ and on
the (non)finiteness of the first few moments. We assume that $X$
is of Pareto-type with index $\alpha>0$ in Cases $1$-$5$ and that
$X$ satisfies $\mu_4<\infty$ in Case $6$. Moreover, the counting
process is supposed to $\mathcal{D}$-average in time to the random
variable $\Lambda$ except for Case~$5$ where it $p$-averages in
time to $\Lambda$.

\vspace{0.1cm}

\paragraph{Case 1: $\alpha\in(0,1)$.}

Since $N(t)\stackrel{a.s.}{\longrightarrow}\infty$ as
$t\rightarrow\infty$, it follows from Theorem~\ref{theo01},
Slutsky's theorem and the continuous mapping theorem that:
\begin{equation*}
  \frac{\widehat{\mathrm{CoVar}(X)}}{\sqrt{N(t)}}=\sqrt{T_{N(t)}-\frac{1}{N(t)}}
  \stackrel{\mathcal{D}}{\longrightarrow}\frac{\sqrt{U_\alpha}}{V_\alpha}
  \quad\text{as}\quad t\rightarrow\infty
\end{equation*}
where the joint distribution of the random vector $(U_\alpha,
V_\alpha)^\prime$ is determined by~(\ref{jointtheo01}).

\vspace{0.1cm}

\paragraph{Case 2: $\alpha=1$, $\mu_1=\infty$.}

Define $(a_t)_{t>0}$ by $\lim_{t\rightarrow\infty}t \,
a_t^{-1}\ell(a_t)=1$ and $(a^\prime_t)_{t>0}$ by
$\lim_{t\rightarrow\infty}t \, a^{\prime
-1}_t\tilde{\ell}(a^\prime_t)=1$ with $\tilde{\ell}(x)=\int_0^x
\frac{\ell(u)}{u} \, du \in \mathrm{RV}_0$. Since
$\frac{a^\prime_t}{a_t}\sim\frac{\left(1/\tilde{\ell}\right)^*(t)}
{\left(1/\ell\right)^*(t)}$ as $t\rightarrow\infty$, where
$\left(\frac{1}{\tilde{\ell}}\right)^*$ and
$\left(\frac{1}{\ell}\right)^*$ respectively stand for the de
Bruyn conjugate of $\frac{1}{\tilde{\ell}}$ and $\frac{1}{\ell}$,
it follows that $\frac{a^\prime_t}{a_t}\in\mathrm{RV}_0$ and
consequently that
$\frac{1}{t}\left(\frac{a^\prime_t}{a_t}\right)^2\rightarrow 0$ as
$t\rightarrow\infty$. Moreover,
$\frac{N(t)}{t}\stackrel{\mathcal{D}}{\longrightarrow}\Lambda$ as
$t\rightarrow\infty$. Hence, Theorem~\ref{theo1} together with
Slutsky's theorem and the continuous mapping theorem gives:
\begin{equation*}
  \frac{a^\prime_t}{a_t} \, \frac{\widehat{\mathrm{CoVar}(X)}}{\sqrt{N(t)}}
  =\sqrt{\left(\frac{a^\prime_t}{a_t}\right)^2 T_{N(t)}-\frac{1}{t}\left(\frac{a^\prime_t}{a_t}\right)^2
  \frac{t}{N(t)}} \stackrel{\mathcal{D}}{\longrightarrow} \frac{\sqrt{U_1}}{\Lambda}
  \quad\text{as}\quad t\rightarrow\infty
\end{equation*}
where the joint distribution of the random vector $(U_1,
\Lambda)^\prime$ is determined by~(\ref{jointtheo1}).

\vspace{0.1cm}

\paragraph{Case 3: $\alpha\in(1,2)$ or $\alpha=1$, $\mu_1<\infty$.}

Define $(a_t)_{t>0}$ by $\lim_{t\rightarrow\infty}t \,
a_t^{-\alpha}\ell(a_t)=1$. Since $\frac{t}{a_t^2}\sim
\frac{a_t^{\alpha-2}}{\ell(a_t)}\rightarrow 0$ and
$\frac{N(t)}{t}\stackrel{\mathcal{D}}{\longrightarrow}\Lambda$ as
$t\rightarrow\infty$, Theorem~\ref{theo12}$(a)$, Slutsky's theorem
and the continuous mapping theorem lead to:
\begin{equation*}
  \frac{\sqrt{N(t)}}{a_t} \, \widehat{\mathrm{CoVar}(X)}=\sqrt{\left(\frac{N(t)}{a_t}\right)^2
  T_{N(t)}-\frac{t}{a_t^2}\frac{N(t)}{t}} \stackrel{\mathcal{D}}{\longrightarrow}
  \frac{1}{\mu_1} \, \sqrt{U_\alpha} \quad\text{as}\quad t\rightarrow\infty
\end{equation*}
where the distribution of $U_\alpha$ is determined by~(\ref{LTUalpha}).\\

Repeating the same arguments as above but using
Theorem~\ref{theo12}$(b)$ instead of Theorem~\ref{theo12}$(a)$, we
also get:
\begin{equation*}
  \frac{t}{a_t} \, \frac{\widehat{\mathrm{CoVar}(X)}}{\sqrt{N(t)}}
  =\sqrt{\left(\frac{t}{a_t}\right)^2 T_{N(t)}-\frac{t}{a_t^2}\frac{t}{N(t)}}
  \stackrel{\mathcal{D}}{\longrightarrow} \frac{1}{\mu_1} \, \frac{\sqrt{U_\alpha}}{\Lambda}
  \quad\text{as}\quad t\rightarrow\infty
\end{equation*}
where the joint distribution of the random vector $(U_\alpha,
\mu_1\Lambda)^\prime$ is determined by~(\ref{jointtheo12}).

\vspace{0.1cm}

\paragraph{Case 4: $\alpha=2$, $\mu_2=\infty$.}

Define $(a^\prime_t)_{t>0}$ by $\lim_{t\rightarrow\infty}t \,
a^{\prime -2}_t \tilde{\ell}(a^\prime_t)=1$ with
$\tilde{\ell}(x)=\int_0^x \frac{\ell(u)}{u} \, du \in
\mathrm{RV}_0$. From $\mu_2=\infty$, it follows that
$\tilde{\ell}(x)\rightarrow\infty$ as $t\rightarrow\infty$ and
consequently that $\frac{t}{a^{\prime
2}_t}\sim\frac{1}{\tilde{\ell}(a^\prime_t)}\rightarrow 0$ as
$t\rightarrow\infty$. Moreover,
$\frac{N(t)}{t}\stackrel{\mathcal{D}}{\longrightarrow}\Lambda$ as
$t\rightarrow\infty$. Thus, Theorem~\ref{theo2}$(a)$, Slutsky's
theorem and the continuous mapping theorem yield:
\begin{equation*}
  \frac{\sqrt{N(t)}}{a^\prime_t} \, \widehat{\mathrm{CoVar}(X)}
  =\sqrt{\left(\frac{N(t)}{a^\prime_t}\right)^2 T_{N(t)}-\frac{t}{a^{\prime 2}_t}\frac{N(t)}{t}}
  \stackrel{\mathcal{D}}{\longrightarrow} \frac{\sqrt{2}}{\mu_1} \, \sqrt{\Lambda}
  \quad\text{as}\quad t\rightarrow\infty.
\end{equation*}

By using Theorem~\ref{theo2}$(b)$ and the arguments above, we also
get:
\begin{equation*}
  \frac{t}{a^\prime_t} \, \frac{\widehat{\mathrm{CoVar}(X)}}{\sqrt{N(t)}}
  =\sqrt{\left(\frac{t}{a^\prime_t}\right)^2 T_{N(t)}-\frac{t}{a^{\prime 2}_t}\frac{t}{N(t)}}
  \stackrel{\mathcal{D}}{\longrightarrow} \frac{\sqrt{2}}{\mu_1} \, \frac{1}{\sqrt{\Lambda}}
  \quad\text{as}\quad t\rightarrow\infty.
\end{equation*}

\paragraph{Case 5: $\alpha\in(2,4)$ or $\alpha=2$, $\mu_2<\infty$.}

Assume that $\{N(t); \, t\geq 0\}$ $p$-averages in time to the
random variable $\Lambda$. From Theorem~\ref{theo24}, we deduce
that $N(t) \, T_{N(t)} \stackrel{p}{\longrightarrow}
\frac{\mu_2}{\mu_1^2}$ as $t\rightarrow\infty$. Using the
continuous mapping theorem, we thus get:
\begin{equation*}
  \widehat{\mathrm{CoVar}(X)}\stackrel{p}{\longrightarrow}\mathrm{CoVar}(X)
  \quad\text{as}\quad t\rightarrow\infty.
\end{equation*}

Moreover, define a sequence $(b_t)_{t>0}$ by
$b_t:=\frac{t^{1-2/\alpha}}{\ell^*(t)}$ where
$\ell^*(t):=\ell^{2/\alpha}\left(\sqrt{c_t}\right)\in\mathrm{RV}_0$
with $(c_t)_{t>0}$ a sequence defined by
$\lim_{t\rightarrow\infty}t \, c_t^{-\alpha/2}
\ell(\sqrt{c_t})=1$. Denote $\sigma^2:=\mathbb{V}X<\infty$ and
consider:
\begin{equation*}
  b_t\left(\widehat{\mathrm{CoVar}(X)}-\mathrm{CoVar}(X)\right)
  =\underbrace{\frac{\mu_1 b_t\left(N(t)\,T_{N(t)}-\frac{\mu_2}{\mu_1^2}\right)}{2\sigma}}_{=:A_t}
  -\underbrace{\frac{\mu_1 b_t\left(N(t)\,T_{N(t)}-\frac{\mu_2}{\mu_1^2}\right)^2}
  {2\sigma\left(\widehat{\mathrm{CoVar}(X)}+\mathrm{CoVar}(X)\right)^2}}_{=:B_t}.
\end{equation*}

From Theorem~\ref{theo24}, we easily deduce by using Slutsky's
theorem that $A_t\stackrel{\mathcal{D}}{\longrightarrow}
\frac{1}{2\mu_1\sigma} \,
\frac{W_{\frac{\alpha}{2}}}{\Lambda^{1-2/\alpha}}$ and that
$B_t\stackrel{p}{\longrightarrow}0$ as $t\rightarrow\infty$
leading by virtue of another application of Slutsky's theorem to:
\begin{equation*}
  \frac{t^{1-2/\alpha}}{\ell^*(t)}\left(\widehat{\mathrm{CoVar}(X)}-\mathrm{CoVar}(X)\right)
  \stackrel{\mathcal{D}}{\longrightarrow} \frac{1}{2\mu_1\sigma} \,
  \frac{W_{\frac{\alpha}{2}}}{\Lambda^{1-2/\alpha}} \quad\text{as}\quad t\rightarrow\infty
\end{equation*}
where $W_{\frac{\alpha}{2}}$ is a stable random variable with
exponent $\alpha/2$ independent of $\Lambda$.

\vspace{0.1cm}

\paragraph{Case 6: $\mu_4<\infty$.} The proof of Theorem~\ref{theo4}
can be repeated using the transformation
$g(x,y)=\sqrt{\frac{y}{x^2}-1}$ and this leads to:
\begin{equation}\label{limmu4C}
  \sqrt{t}\left(\widehat{\mathrm{CoVar}(X)}-\mathrm{CoVar}(X)\right)
  \stackrel{\mathcal{D}}{\longrightarrow}\frac{\mathrm{N}\left(0,\frac{\sigma_*^2\mu_1^2}{4\sigma^2}\right)}
  {\sqrt{\Lambda}}\quad\text{as}\quad t\rightarrow\infty
\end{equation}
where the random variable
$\mathrm{N}\left(0,\frac{\sigma_*^2\mu_1^2}{4\sigma^2}\right)$ is
independent of $\Lambda$ and has the normal distribution with mean
$0$ and variance $\frac{\sigma_*^2\mu_1^2}{4\sigma^2}$, with
$\sigma_*^2$ defined by~(\ref{defsigma*}) and
$\sigma^2:=\mathbb{V}X<\infty$.\\

Assume $\mathbb{E}\left\{\Lambda^{-1}\right\}<\infty$. When
$t\left(\widehat{\mathrm{CoVar}(X)}-\mathrm{CoVar}(X)\right)^2$ is
uniformly integrable, the first two moments of the limiting
distribution in~(\ref{limmu4C}) permit to determine the limiting
behavior of
$\mathrm{CoVar}\left(\widehat{\mathrm{CoVar}(X)}\right)$. Indeed,
on the one hand:
\begin{equation*}
  \mathbb{E}\left\{\sqrt{t}\left(\widehat{\mathrm{CoVar}(X)}-\mathrm{CoVar}(X)\right)\right\}
  =\sqrt{t}\left(\mathbb{E}\left\{\widehat{\mathrm{CoVar}(X)}\right\}-\mathrm{CoVar}(X)\right)
  \rightarrow 0 \quad\text{as}\quad t\rightarrow\infty
\end{equation*}
which leads to:
\begin{equation}\label{AUEcov}
  \mathbb{E}\left\{\widehat{\mathrm{CoVar}(X)}\right\}\rightarrow\mathrm{CoVar}(X)
  \quad\text{as}\quad t\rightarrow\infty.
\end{equation}

One the other hand, we also get:
\begin{equation*}
  \mathbb{V}\left\{\sqrt{t}\left(\widehat{\mathrm{CoVar}(X)}-\mathrm{CoVar}(X)\right)\right\}
  =t \, \mathbb{V} \left\{\widehat{\mathrm{CoVar}(X)}\right\}
  \sim \frac{\sigma_*^2\mu_1^2 \, \mathbb{E}\left\{\Lambda^{-1}\right\}}{4\sigma^2}
  \quad\text{as}\quad t\rightarrow\infty
\end{equation*}
implying that:
\begin{equation*}
  \mathbb{V}\left\{\widehat{\mathrm{CoVar}(X)}\right\}
  \sim \frac{\sigma_*^2\mu_1^2\,\mathbb{E}\left\{\Lambda^{-1}\right\}}{4\sigma^2}\,\frac{1}{t}
  \quad\text{as}\quad t\rightarrow\infty.
\end{equation*}

Consequently, under the above uniform integrability condition, the
coefficient of variation of the sample coefficient of variation
asymptotically behaves as:
\begin{equation*}
  \mathrm{CoVar}\left(\widehat{\mathrm{CoVar}(X)}\right)
  \sim \frac{\sigma_*\mu_1^2\,\sqrt{\mathbb{E}\left\{\Lambda^{-1}\right\}}}{2\sigma^2}\,\frac{1}{\sqrt{t}}
  \quad\text{as}\quad t\rightarrow\infty.
\end{equation*}

In addition, it results from~(\ref{limmu4C}) and~(\ref{AUEcov})
that $\widehat{\mathrm{CoVar}(X)}$ is a consistent and
asymptotically unbiased estimator for $\mathrm{CoVar}(X)$.

\subsection{Sample Dispersion}\label{sec_appli_2}

Adapting the results of Section~\ref{sec_weaklaws} to the random
variable $C_{N(t)}$ defined in~(\ref{defC}) permits us to derive
limiting distributions for the appropriately normalized sample of
dispersion $\widehat{\mathrm{D}(X)}$ thanks to~(\ref{reldisp}).
Different cases are considered as for the sample coefficient of
variation. We assume that $X$ is of Pareto-type with index
$\alpha>0$ in Cases $1$-$5$ and that $X$ satisfies $\mu_4<\infty$
in Case $6$. Moreover, the counting process is supposed to
$\mathcal{D}$-average in time to the random variable $\Lambda$
except for Case~$5$ where it $p$-averages in time to $\Lambda$.

\vspace{0.1cm}

\paragraph{Case 1: $\alpha\in(0,1)$.}

Define $(a_t)_{t>0}$ by $\lim_{t\rightarrow\infty}t \,
a_t^{-\alpha}\ell(a_t)=1$. It follows from Theorem~\ref{theo01}
and the continuous mapping theorem that:
\begin{equation*}
  \frac{1}{a_t} \, C_{N(t)} \stackrel{\mathcal{D}}{\longrightarrow}
  \frac{U_\alpha}{V_\alpha} \quad\text{as}\quad t\rightarrow\infty
\end{equation*}
where the joint distribution of the random vector $(U_\alpha,
V_\alpha)^\prime$ is determined by~(\ref{jointtheo01}).\\

Since $N(t)\stackrel{a.s.}{\longrightarrow}\infty$ and
$\frac{1}{a_t}\sum_{i=1}^{N(t)}X_i\stackrel{\mathcal{D}}{\longrightarrow}V_\alpha$
as $t\rightarrow\infty$, where the distribution of $V_\alpha$ is
determined by~(\ref{LTValpha}), Slutsky's theorem together with
the continuous mapping theorem yields:
\begin{equation*}
  \frac{1}{a_t} \, \widehat{\mathrm{D}(X)}
  =\frac{1}{a_t} \, C_{N(t)}-\frac{1}{N(t)} \frac{1}{a_t}\sum_{i=1}^{N(t)}X_i
  \stackrel{\mathcal{D}}{\longrightarrow} \frac{U_\alpha}{V_\alpha}
  \quad\text{as}\quad t\rightarrow\infty.
\end{equation*}

\paragraph{Case 2: $\alpha=1$, $\mu_1=\infty$.}

Define $(a_t)_{t>0}$ by $\lim_{t\rightarrow\infty}t \, a_t^{-1}
\ell(a_t)=1$ and $(a^\prime_t)_{t>0}$ by
$\lim_{t\rightarrow\infty}t \, a^{\prime -1}_t
\tilde{\ell}(a^\prime_t)=1$ with $\tilde{\ell}(x)=\int_0^x
\frac{\ell(u)}{u} \, du \in \mathrm{RV}_0$. It follows from
Theorem~\ref{theo1} and the continuous mapping theorem that:
\begin{equation*}
  \frac{a^\prime_t}{a_t^2} \, C_{N(t)} \stackrel{\mathcal{D}}{\longrightarrow}
  \frac{U_1}{\Lambda} \quad\text{as}\quad t\rightarrow\infty
\end{equation*}
where the joint distribution of the random vector $(U_1,
\Lambda)^\prime$ is determined by~(\ref{jointtheo1}).\\

Since
$\frac{a^\prime_t}{a_t}\sim\frac{\left(1/\tilde{\ell}\right)^*(t)}
{\left(1/\ell\right)^*(t)}$ as $t\rightarrow\infty$, where
$\left(\frac{1}{\tilde{\ell}}\right)^*$ and
$\left(\frac{1}{\ell}\right)^*$ respectively stand for the de
Bruyn conjugate of $\frac{1}{\tilde{\ell}}$ and $\frac{1}{\ell}$,
it follows that $\frac{a^\prime_t}{a_t}\in\mathrm{RV}_0$ and
consequently that
$\frac{1}{t}\left(\frac{a^\prime_t}{a_t}\right)^2\rightarrow 0$ as
$t\rightarrow\infty$. Moreover, using the same independence and
conditioning arguments as in the proof of Theorem~\ref{theo24}, we
obtain that at any continuity points $x$ and $y$ of the
distribution function of $\Lambda$:
\begin{equation*}
  \mathbb{P}\left[\frac{N(t)}{t}\leq x, \frac{1}{a^\prime_t}\sum_{i=1}^{N(t)}X_i\leq y\right]
  \rightarrow\mathbb{P}\left[\Lambda\leq x\right]\,\mathbb{P}\left[\Lambda\leq y\right]
  \quad\text{as}\quad t\rightarrow\infty
\end{equation*}
i.e., since
$\frac{N(t)}{t}\stackrel{\mathcal{D}}{\longrightarrow}\Lambda$ and
$\frac{1}{a^\prime_t}\sum_{i=1}^{N(t)}X_i\stackrel{\mathcal{D}}{\longrightarrow}\Lambda$
as $t\rightarrow\infty$, that:
\begin{equation*}
  \left(\frac{N(t)}{t}, \, \frac{1}{a^\prime_t} \sum_{i=1}^{N(t)}X_i\right)^\prime
  \stackrel{\mathcal{D}}{\longrightarrow} (\Lambda, \Lambda^*)^\prime \quad\text{as}\quad t\rightarrow\infty
\end{equation*}
where $\Lambda^*$ is an independent copy of $\Lambda$.\\

Using the continuous mapping theorem, we thus deduce:
\begin{equation*}
  \frac{t}{N(t)}\frac{\sum_{i=1}^{N(t)}X_i}{a^\prime_t}
  \stackrel{\mathcal{D}}{\longrightarrow} \frac{\Lambda^*}{\Lambda} \quad\text{as}\quad t\rightarrow\infty.
\end{equation*}

Hence, Slutsky's theorem gives:
\begin{equation*}
  \frac{a^\prime_t}{a_t^2} \, \widehat{\mathrm{D}(X)}=\frac{a^\prime_t}{a_t^2} \, C_{N(t)}-
  \frac{1}{t}\left(\frac{a^\prime_t}{a_t}\right)^2\frac{t}{N(t)}\frac{\sum_{i=1}^{N(t)}X_i}{a^\prime_t}
  \stackrel{\mathcal{D}}{\longrightarrow} \frac{U_1}{\Lambda}
  \quad\text{as}\quad t\rightarrow\infty.
\end{equation*}

\paragraph{Case 3: $\alpha\in(1,2)$ or $\alpha=1$, $\mu_1<\infty$.}

Define $(a_t)_{t>0}$ by $\lim_{t\rightarrow\infty}t \,
a_t^{-\alpha}\ell(a_t)=1$. It follows from
Theorem~\ref{theo12}$(a)$ and the continuous mapping theorem that:
\begin{equation*}
  \frac{N(t)}{a_t^2} \, C_{N(t)} \stackrel{\mathcal{D}}{\longrightarrow}
  \frac{1}{\mu_1} \, U_\alpha \quad\text{as}\quad t\rightarrow\infty
\end{equation*}
where the distribution of $U_\alpha$ is determined by~(\ref{LTUalpha}).\\

Since $\frac{t}{a_t^2}\sim
\frac{a_t^{\alpha-2}}{\ell(a_t)}\rightarrow 0$,
$\frac{N(t)}{t}\stackrel{\mathcal{D}}{\longrightarrow}\Lambda$ and
$\overline{X}\stackrel{p}{\longrightarrow}\mu_1$ as
$t\rightarrow\infty$, Slutsky's theorem leads to:
\begin{equation*}
  \frac{N(t)}{a_t^2} \, \widehat{\mathrm{D}(X)}
  =\frac{N(t)}{a_t^2} \, C_{N(t)}- \frac{t}{a_t^2}\frac{N(t)}{t} \,\overline{X}
  \stackrel{\mathcal{D}}{\longrightarrow}
  \frac{1}{\mu_1} \, U_\alpha \quad\text{as}\quad t\rightarrow\infty.
\end{equation*}

Repeating the same arguments as above but using
Theorem~\ref{theo12}$(b)$ instead of Theorem~\ref{theo12}$(a)$, we
also get:
\begin{equation*}
  \frac{t}{a_t^2} \, \widehat{\mathrm{D}(X)}=
  \frac{t}{a_t^2} \, C_{N(t)} - \frac{t}{a_t^2} \, \overline{X}
  \stackrel{\mathcal{D}}{\longrightarrow} \frac{1}{\mu_1} \, \frac{U_\alpha}{\Lambda}
  \quad\text{as}\quad t\rightarrow\infty
\end{equation*}
where the joint distribution of the random vector $(U_\alpha,
\mu_1\Lambda)^\prime$ is determined by~(\ref{jointtheo12}).

\vspace{0.1cm}

\paragraph{Case 4: $\alpha=2$, $\mu_2=\infty$.}

Define $(a^\prime_t)_{t>0}$ by $\lim_{t\rightarrow\infty}t \,
a^{\prime -2}_t \tilde{\ell}(a^\prime_t)=1$ with
$\tilde{\ell}(x)=\int_0^x \frac{\ell(u)}{u} \, du \in
\mathrm{RV}_0$. It follows from Theorem~\ref{theo2}$(a)$ and the
continuous mapping theorem that:
\begin{equation*}
  \frac{N(t)}{a^{\prime 2}_t} \, C_{N(t)} \stackrel{\mathcal{D}}{\longrightarrow}
  \frac{2}{\mu_1} \, \Lambda \quad\text{as}\quad t\rightarrow\infty.
\end{equation*}

From $\mu_2=\infty$, it follows that
$\tilde{\ell}(x)\rightarrow\infty$ as $t\rightarrow\infty$ and
consequently that $\frac{t}{a^{\prime
2}_t}\sim\frac{1}{\tilde{\ell}(a^\prime_t)}\rightarrow 0$ as
$t\rightarrow\infty$. Moreover, since
$\frac{N(t)}{t}\stackrel{\mathcal{D}}{\longrightarrow}\Lambda$ and
$\overline{X}\stackrel{p}{\longrightarrow}\mu_1$ as
$t\rightarrow\infty$, Slutsky's theorem yields:
\begin{equation*}
  \frac{N(t)}{a^{\prime 2}_t} \, \widehat{\mathrm{D}(X)}
  =\frac{N(t)}{a^{\prime 2}_t} \, C_{N(t)} - \frac{t}{a^{\prime 2}_t} \frac{N(t)}{t} \, \overline{X}
  \stackrel{\mathcal{D}}{\longrightarrow} \frac{2}{\mu_1} \, \Lambda
  \quad\text{as}\quad t\rightarrow\infty.
\end{equation*}

By using Theorem~\ref{theo2}$(b)$ and the arguments above, we also
get:
\begin{equation*}
  \frac{t}{a^{\prime 2}_t} \, \widehat{\mathrm{D}(X)}
  =\frac{t}{a^{\prime 2}_t} \, C_{N(t)} - \frac{t}{a^{\prime 2}_t} \, \overline{X}
  \stackrel{\mathcal{D}}{\longrightarrow} \frac{2}{\mu_1} \quad\text{as}\quad t\rightarrow\infty.
\end{equation*}

\paragraph{Case 5: $\alpha\in(2,4)$ or $\alpha=2$, $\mu_2<\infty$.}

Assume that $\{N(t); \, t\geq 0\}$ $p$-averages in time to the
random variable $\Lambda$. Define a sequence $(b_t)_{t>0}$ by
$b_t:=\frac{t^{1-2/\alpha}}{\ell^*(t)}$ where
$\ell^*(t):=\ell^{2/\alpha}\left(\sqrt{c_t}\right)\in\mathrm{RV}_0$
with $(c_t)_{t>0}$ a sequence defined by
$\lim_{t\rightarrow\infty}t \, c_t^{-\alpha/2} \ell(\sqrt{c_t})=1$
and consider:
\begin{eqnarray*}
  b_t\left(\widehat{\mathrm{D}(X)}-\mathrm{D}(X)\right)
  &=& b_t\left(C_{N(t)}-\frac{\mu_2}{\mu_1}-\overline{X}+\mu_1\right) \cr
  &=&
  \underbrace{\frac{b_t}{\overline{X}}\left(\frac{1}{N(t)}\sum_{i=1}^{N(t)} X_i^2 - \mu_2\right)}_{=:A_t}
  -\underbrace{b_t\left(\overline{X}-\mu_1\right)\left(1+\frac{\mu_2}{\mu_1}
  \frac{1}{\overline{X}}\right)}_{=:B_t}.
\end{eqnarray*}

Since $\overline{X}\stackrel{p}{\longrightarrow}\mu_1$,
$\frac{N(t)}{t}\stackrel{p}{\longrightarrow}\Lambda$ and
$\frac{\ell^*(N(t))}{\ell^*(t)}\stackrel{p}{\longrightarrow}1$ as
$t\rightarrow\infty$, Slutsky's theorem and the continuous mapping
theorem give:
\begin{equation*}
  A_t=\frac{1}{\overline{X}}\left(\frac{t}{N(t)}\right)^{1-2/\alpha}\frac{\ell^*(N(t))}{\ell^*(t)}
  \frac{N(t)^{1-2/\alpha}}{\ell^*(N(t))}
  \left(\frac{1}{N(t)}\sum_{i=1}^{N(t)} X_i^2 - \mu_2\right)
  \stackrel{\mathcal{D}}{\longrightarrow} \frac{1}{\mu_1}\frac{W_{\frac{\alpha}{2}}}{\Lambda^{1-2/\alpha}}
  \quad\text{as}\quad t\rightarrow\infty
\end{equation*}
thanks to~(\ref{AnsCLTW}) and the independence of $\Lambda$ and
$W_{\frac{\alpha}{2}}$ which is a stable random variable with
exponent $\alpha/2$.\\

Incidentally, using~(\ref{AnsCLTN}) yields:
\begin{equation*}
  \sqrt{t}\left(\frac{1}{N(t)}\sum_{i=1}^{N(t)}X_i - \mu_1\right)
  =\sqrt{\frac{t}{N(t)}}\sqrt{N(t)}\left(\frac{1}{N(t)}\sum_{i=1}^{N(t)}X_i - \mu_1\right)
  \stackrel{\mathcal{D}}{\longrightarrow} \frac{\mathrm{N}(0,\sigma^2)}{\sqrt{\Lambda}}
  \quad\text{as}\quad t\rightarrow\infty
\end{equation*}
since the random variable $\mathrm{N}(0,\sigma^2)$ which has the
normal distribution with mean $0$ and variance
$\sigma^2:=\mathbb{V}X<\infty$ is independent of $\Lambda$.
Consequently, by considering the equality:
\begin{equation*}
  b_t\left(\overline{X}-\mu_1\right)=\frac{t^{1/2-2/\alpha}}{\ell^*(t)}
  \, \sqrt{t}\left(\frac{1}{N(t)}\sum_{i=1}^{N(t)} X_i-\mu_1\right)
\end{equation*}
and since $\frac{t^{1/2-2/\alpha}}{\ell^*(t)}\rightarrow 0$ and
$\overline{X}\stackrel{p}{\longrightarrow}\mu_1$ as
$t\rightarrow\infty$, Slutsky's theorem together with the
continuous mapping theorem implies that
$B_t\stackrel{p}{\longrightarrow}0$ as $t\rightarrow\infty$. By
virtue of another application of Slutsky's theorem, we hence
obtain:
\begin{equation*}
  \frac{t^{1-2/\alpha}}{\ell^*(t)}\left(\widehat{\mathrm{D}(X)}-\mathrm{D}(X)\right)
  \stackrel{\mathcal{D}}{\longrightarrow} \frac{1}{\mu_1} \,
  \frac{W_{\frac{\alpha}{2}}}{\Lambda^{1-2/\alpha}} \quad\text{as}\quad t\rightarrow\infty.
\end{equation*}

The latter relation shows in particular that:
\begin{equation*}
  \widehat{\mathrm{D}(X)}\stackrel{p}{\longrightarrow}\mathrm{D}(X)\quad\text{as}\quad t\rightarrow\infty.
\end{equation*}

\paragraph{Case 6: $\mu_4<\infty$.} Using $g(x,y)=\frac{y}{x}-x$ in the proof of
Theorem~\ref{theo4} yields:
\begin{equation}\label{limmu4D}
  \sqrt{t}\left(\widehat{\mathrm{D}(X)}-\mathrm{D}(X)\right)
  \stackrel{\mathcal{D}}{\longrightarrow} \frac{\mathrm{N}\left(0,\sigma_{**}^2\right)}
  {\sqrt{\Lambda}} \quad\text{as}\quad t\rightarrow\infty
\end{equation}
where the random variable $\mathrm{N}\left(0,\sigma_{**}^2\right)$
is independent of $\Lambda$ and has the normal distribution with
mean $0$ and variance $\sigma_{**}^2$ defined by:
\begin{equation*}
  \sigma_{**}^2:=\mu_2-\mu_1^2+\frac{\mu_2^3}{\mu_1^4}-2\frac{\mu_3}{\mu_1}
  -2\frac{\mu_2\mu_3}{\mu_1^3}+2\left(\frac{\mu_2}{\mu_1}\right)^2+\frac{\mu_4}{\mu_1^2}.
\end{equation*}

Assume $\mathbb{E}\left\{\Lambda^{-1}\right\}<\infty$. When
$t\left(\widehat{\mathrm{D}(X)}-\mathrm{D}(X)\right)^2$ is
uniformly integrable, the first two moments of the limiting
distribution in~(\ref{limmu4D}) permit to determine the limiting
behavior of $\mathrm{D}\left(\widehat{\mathrm{D}(X)}\right)$.
Indeed, on the one hand:
\begin{equation*}
  \mathbb{E}\left\{\sqrt{t}\left(\widehat{\mathrm{D}(X)}-\mathrm{D}(X)\right)\right\}
  =\sqrt{t}\left(\mathbb{E}\left\{\widehat{\mathrm{D}(X)}\right\}-\mathrm{D}(X)\right)
  \rightarrow 0 \quad\text{as}\quad t\rightarrow\infty
\end{equation*}
leading to:
\begin{equation}\label{AUEdisp}
  \mathbb{E}\left\{\widehat{\mathrm{D}(X)}\right\}\rightarrow\mathrm{D}(X)
  \quad\text{as}\quad t\rightarrow\infty.
\end{equation}

Note that~(\ref{limmu4D}) together with~(\ref{AUEdisp}) implies
that $\widehat{\mathrm{D}(X)}$ is a consistent and asymptotically
unbiased estimator for $\mathrm{D}(X)$.\\

On the other hand, we also get:
\begin{equation*}
  \mathbb{V}\left\{\sqrt{t}\left(\widehat{\mathrm{D}(X)}-\mathrm{D}(X)\right)\right\}
  =t \, \mathbb{V} \left\{\widehat{\mathrm{D}(X)}\right\}
  \sim \sigma_{**}^2 \, \mathbb{E}\left\{\Lambda^{-1}\right\}
  \quad\text{as}\quad t\rightarrow\infty
\end{equation*}
which implies:
\begin{equation*}
  \mathbb{V}\left\{\widehat{\mathrm{D}(X)}\right\}
  \sim \sigma_{**}^2 \, \mathbb{E}\left\{\Lambda^{-1}\right\} \, \frac{1}{t}
  \quad\text{as}\quad t\rightarrow\infty.
\end{equation*}

Consequently, under the above uniform integrability condition, the
dispersion of the sample dispersion asymptotically behaves as:
\begin{equation*}
  \mathrm{D}\left(\widehat{\mathrm{D}(X)}\right)
  \sim \frac{\sigma^2_{**} \, \mu_1 \, \mathbb{E}\left\{\Lambda^{-1}\right\}}{\sigma^2} \, \frac{1}{t}
  \quad\text{as}\quad t\rightarrow\infty
\end{equation*}
where $\sigma^2:=\mathbb{V}X<\infty$.

\section{Conclusion}\label{ccl}

We have derived limits in distribution for the random variable
$T_{N(t)}$ defined in~(\ref{defT}) when the distribution function
$F$ of $X_1$ is of Pareto-type with index $\alpha>0$ or is such
that $\mu_4<\infty$. Furthermore, the counting process $\{N(t); \,
t\geq 0\}$ has been chosen to $\mathcal{D}$-average in time or, in
a single case, to $p$-average in time. Different results have
shown up according to the range of $\alpha$ and to the
(non)finiteness of the first moments. These results have then been
used to analyze limiting properties of two risk measures that are
used in many applications, namely the sample coefficient of
variation and the sample dispersion. However, the practical
applicability of these risk measures depends on the existence of
sufficiently many moments of the underlying distribution. Hence,
the results we have obtained help at illustrating their asymptotic
behavior when such a moment condition is not satisfied.\\

We point out that by choosing the limiting random variable
$\Lambda$ to be degenerate at the point $1$ in our results, we
retrieve results of Albrecher and Teugels~\cite{at04} where the
counting process is assumed to be deterministic.


\begin{thebibliography}{10}

\bibitem{at04}
H.~Albrecher and J.L. Teugels.
\newblock Asymptotic analysis of measures of variation.
\newblock Technical Report 2004-042, EURANDOM, Technical University of
  Eindhoven, The Netherlands, 2004.

\bibitem{bgst04}
J.~Beirlant, Y.~Goegebeur, J.~Segers, and J.L. Teugels.
\newblock {\em {S}tatistics of {E}xtremes: {T}heory and {A}pplications}.
\newblock John Wiley \& Sons, Chichester, 2004.

\bibitem{Bill68}
P.~Billingsley.
\newblock {\em {C}onvergence of {P}robability {M}easures}.
\newblock John Wiley \& Sons, New York, 1968.

\bibitem{BGT}
N.H. Bingham, C.M. Goldie, and J.L. Teugels.
\newblock {\em {R}egular {V}ariation}.
\newblock Cambridge University Press, Cambridge, 1987.

\bibitem{Chung01}
K.L. Chung.
\newblock {\em A {C}ourse in {P}robability {T}heory}.
\newblock Third edition. Academic Press, San Diego, 2001.

\bibitem{Feller2}
W.~Feller.
\newblock {\em {A}n {I}ntroduction to {P}robability {T}heory and {I}ts
  {A}pplications. {V}ol. {II}}.
\newblock Second edition. John Wiley \& Sons, New York, 1971.

\bibitem{gr97}
J.~Grandell.
\newblock {\em Mixed {P}oisson {P}rocesses}.
\newblock Monographs on Statistics and Applied Probability 77. Chapman \& Hall,
  London, 1997.

\bibitem{L05}
S.A. Ladoucette.
\newblock Asymptotic results on the moments of the ratio of the random sum of
  squares to the square of the random sum.
\newblock Technical Report 2005-03, University Center for Statistics,
  Department of Mathematics, Catholic University of Leuven, Belgium, 2005.

\bibitem{MAC97}
T.~Mack.
\newblock Schadenversicherungsmathematik.
\newblock {\em Verlag Versicherungswirtschaft e.V., Karlsruhe}, 1997.

\bibitem{SERF}
R.J. Serfling.
\newblock {\em Approximation {T}heorems of {M}athematical {S}tatistics}.
\newblock John Wiley \& Sons, New York, 1980.

\end{thebibliography}
\end{document}